# La controverse de 1874 entre Camille Jordan et Leopold Kronecker.

**Frédéric Brechenmacher (*).**

**Résumé.** Une vive querelle oppose en 1874 Camille Jordan et Leopold Kronecker sur l'organisation de la théorie des formes bilinéaires, considérée comme permettant un traitement « général » et « homogène » de nombreuses questions développées dans des cadres théoriques variés au XIX$^{e}$ siècle et dont le problème principal est reconnu comme susceptible d'être résolu par deux théorèmes énoncés indépendamment par Jordan et Weierstrass. Cette controverse, suscitée par la rencontre de deux théorèmes que nous considèrerions aujourd'hui équivalents, nous permettra de questionner l'identité algébrique de pratiques polynomiales de manipulations de « formes » mises en œuvre sur une période antérieure aux approches structurelles de l'algèbre linéaire qui donneront à ces pratiques l'identité de méthodes de caractérisation des classes de similitudes de matrices. Nous montrerons que les pratiques de réductions canoniques et de calculs d'invariants opposées par Jordan et Kronecker manifestent des identités multiples indissociables d'un contexte social daté et qui dévoilent des savoirs tacites, des modes de pensées locaux mais aussi, au travers de regards portés sur une histoire à long terme impliquant des travaux d'auteurs comme Lagrange, Laplace, Cauchy ou Hermite, deux philosophies internes sur la signification de la généralité indissociables d'idéaux disciplinaires opposant algèbre et arithmétique. En questionnant les identités culturelles de telles pratiques cet article vise à enrichir l'histoire de l'algèbre linéaire, souvent abordée dans le cadre de problématiques liées à l'émergence de structures et par l'intermédiaire de l'histoire d'une théorie, d'une notion ou d'un mode de raisonnement.

**Abstract.** Throughout the whole year of 1874, Camille Jordan and Leopold Kronecker were quarrelling over the organisation of the theory of bilinear forms, a theory that was considered as giving a new "homogeneous" and "general" treatment to different problems referring to various theories developed in the 19$^{th}$ century. Two theorems stated independently by Jordan and Weierstrass could be used to solve the main problem of the theory and although these theorems would be considered equivalent as regard to modern mathematics it was their opposition that generated the 1874 controversy. As we will be looking into this quarrel, our purpose will be to discuss the algebraic identities of practices used before the time of linear algebra when these practices would be seen as methods for the classification of similar matrices. Studying the complex identities of practices such as Jordan's canonical reduction and Kronecker's invariant computation sheds some light on some cultural and context issues such as tacit knowledge or local ways of thinking and therefore aims at a deeper understanding of the history of linear algebra without focusing on issues related to the origins of theories or structures. The different ways the two opponents referred to the history of a hundred years old mechanical problem and the works of Lagrange, Laplace, Cauchy and Hermite highlight an opposition of two disciplinary ideals on algebra and arithmetic as well as two internal philosophies of "generality".

---

* F. BRECHENMACHER. Laboratoire de Mathématiques Lens (LML, EA2462). Fédération de Recherche Mathématique du Nord-Pas-de-Calais (CNRS, FR 2956). Université d'Artois (IUFM du Nord Pas de Calais).
Faculté des Sciences Jean Perrin, rue Jean Souvraz S.P. 18, 62 307 Lens Cedex France.
Courrier électronique : frederic.brechenmacher@euler.univ-artois.fr.

Trois étapes permettent de décrire l'évolution de la controverse qui oppose, en 1874, Camille Jordan et Leopold Kronecker. Deux communications successives adressées aux académies de Paris et Berlin marquent l'origine d'une querelle de priorité qui suscite, dans la sphère privée, un échange épistolaire durant l'hiver 1874 puis, au printemps de cette même année, une série de notes et mémoires qui sont autant d'attaques et contre attaques au sein de ce que Jordan désigne comme la « scène publique ». En suivant cette évolution chronologique, détaillée en annexe 1 et qui nous amènera à remonter dans le temps lorsque Jordan et Kronecker mobiliseront des références à des textes anciens, nous verrons que la question de priorité posée par la querelle oppose deux théorèmes. L'un, énoncé par Jordan dans le cadre de son *Traité des substitutions et des équations algébriques* de 1870, établit une forme canonique des substitutions du groupe linéaire (annexe 2). L'autre, établi par Karl Weierstrass en 1868, énonce un système d'invariants, les diviseurs élémentaires, des couples non singuliers de formes bilinéaires (annexe 3). Si la question de l'identité entre ces deux théorèmes est un élément essentiel de la controverse ([1]), les arguments qu'opposent les deux protagonistes ne se limitent pas à une question d'équivalence mathématique et permettent de porter un regard sur l'histoire de l'algèbre linéaire qui dévoile des pratiques algébriques antérieures aux approches structurelles comme la théorie des matrices des années trente du XX^e^ siècle. Notre problématique vise à poser l'identité algébrique de telles pratiques comme « un problème et non une tautologie » pour reprendre l'expression employée par Catherine Goldstein qui, dans son ouvrage intitulé *Un théorème de Fermat et ses lecteurs,* a montré la pertinence de la question d'identité pour décrire des évolutions qui ne relèvent pas simplement d'une activité de recherche de nouveaux résultats ou de meilleures preuves mais témoignent « de pratiques et, à travers elles, de la manière dont est estimée l'innovation » [Goldstein 1995, p. 16] . Distinguer entre une « version moderne » et une « version originale » des deux théorèmes opposés en 1874 pose une question d'identité semblable à celle développée par Hourya Sinaceur dans son histoire du théorème de Sturm et qui nous amènera à « revenir, par delà les traditions didactiques, aux mémoires originaux. On y apprend toutes les identités que le "progrès" efface : identité d'un contexte, d'un objectif, d'une perspective, d'un langage, sans parler de tout ce qui reste implicite sans manquer d'être là » [Sinaceur 1991, p. 21]. Nous verrons que la richesse de l'histoire d'un théorème d'algèbre linéaire comme le théorème de Jordan de la décomposition matricielle provient souvent de ce qui échappe aux descriptions exprimées dans le cadre de mathématiques qui nous sont contemporaines, C'est pourquoi nous renvoyons systématiquement nos commentaires modernes aux notes de bas de page. ([2]).

---

[1] D'un point de vue qui nous est contemporain, le théorème de réduction d'une matrice à coefficients complexes à sa forme canonique de Jordan est équivalent au théorème des diviseurs élémentaires. Consulter par exemple le manuel de [Gantmacher 1959] ainsi que les annexes 2 et 3. Ci-dessous, trois exemples de décompositions matricielles associées à la décomposition en diviseurs élémentaires d'un même polynôme caractéristique $|A-\lambda I|=(\lambda-1)^2(\lambda-2)^3(\lambda-3)$.

| | | | |
|---|---|---|---|
| Formes de Jordan. | $\begin{pmatrix} 1 &&&&& \\ & 1 &&&& \\ && 2 &&& \\ &&& 2 && \\ &&&& 2 & \\ &&&&& 3 \end{pmatrix}$ | $\begin{pmatrix} 1 &&&&& \\ & 1 &&&& \\ && 2 & 1 && \\ &&& 2 & 1 & \\ &&&& 2 & \\ &&&&& 3 \end{pmatrix}$ | $\begin{pmatrix} 1 &&&&& \\ & 1 &&&& \\ && 2 & 0 && \\ &&& 2 & 1 & \\ &&&& 2 & \\ &&&&& 3 \end{pmatrix}$ |
| Diviseurs élémentaires. | $(\lambda-1), (\lambda-1),$ $(\lambda-2), (\lambda-2), (\lambda-2),$ $(\lambda-3)$ | $(\lambda-1), (\lambda-1),$ $(\lambda-2)^3,$ $(\lambda-3)$ | $(\lambda-1), (\lambda-1),$ $(\lambda-2), (\lambda-2)^2,$ $(\lambda-3)$ |

[2] En complément du cas général de l'équivalence des couples de matrices considéré ici, une étude détaillée du cas symétrique est proposée dans une publication ([Brechenmacher 200?]) s'appuyant sur la querelle entre Jordan et Kronecker comme moment de référence pour questionner, sur le long terme (1766-1874), l'identité algébrique d'une pratique propre à ce que nous désignerons plus loin comme la *discussion sur l'équation des petites oscillations*. Afin que ces deux publications, dont l'intersection est non vide, s'articulent tout en restant chacune

# I. Origines d'une controverse.

## 1. Querelle de priorité et organisation de la théorie des formes bilinéaires.

La controverse a pour origine l'ambition de Jordan, formulée dans une note aux *Comptes rendus* du 22 décembre 1873, de réorganiser la théorie des « polynômes bilinéaires » par des méthodes de réduction à des formes canoniques. S'intercalant entre la note de Jordan et la parution dans le *Journal de Liouville* du mémoire « Sur les formes bilinéaires » qu'elle annonce ([3]), Kronecker réplique le 19 janvier 1874 par une lecture à l'Académie de Berlin d'un mémoire intitulé « Uber Schaaren von quadratischen und bilinearen Formen ». Comme nous allons le détailler dans ce paragraphe, les deux géomètres campent, dès ces premières communications, des positions sur une querelle de priorité indissociable d'un enjeu important qui est d'organiser l'objet et les méthodes de la théorie des formes bilinéaires alors que des applications récentes à la géométrie, l'arithmétique des formes quadratiques et divers problèmes d'intégrations de systèmes différentiels annoncent le rôle essentiel que jouera cette théorie dans les mathématiques de la fin du XIX^e^ siècle ([4]). Il est au préalable nécessaire de dresser un panorama des principales étapes ayant jalonné le développement de la théorie des formes bilinéaires par un petit groupe de géomètres berlinois dans les années 1860 ainsi que des éléments ayant amené Jordan à élaborer une réorganisation de cette théorie ([5]). Deux mémoires successifs, publiés en 1866 par Elwin Christoffel et Kronecker dans le *Journal de Crelle*, avaient jeté les bases d'une « théorie » portant sur la caractérisation des formes bilinéaires – étant données deux formes $A = \Sigma A_{\alpha\beta} x_\alpha y_\beta$ et $A' = \Sigma A'_{\alpha\beta} x_\alpha y_\beta$ , déterminer s'il est possible de transformer l'une en l'autre par des substitutions linéaires opérées sur les deux systèmes de variables – par des méthodes de calculs d'invariants. Dans le cadre de ses travaux sur les fonctions elliptiques et abéliennes, et plus particulièrement sur les transformations des fonctions thêta de plusieurs variables, Kronecker avait été amené à considérer le problème de la transformation simultanée de deux formes $A$ et $B$. La recherche des invariants d'un couple $(A, B)$ était abordée par autonome, les résultats de l'une ont été résumés dans l'autre. Ces travaux sont issus d'une thèse de doctorat, menée sous la direction de J. Dhombres à l'EHESS et au centre A. Koyré, posant plus généralement la question des différentes identités du théorème de Jordan sur la période 1870-1930 [Brechenmacher 2006a]. Au sujet du théorème des diviseurs élémentaires, consulter les travaux menés par Thomas Hawkins sur l'histoire de la théorie des matrices et qui ont été à la base de nos propres recherches [Hawkins, 1977]. D'autres sources d'inspirations ont été les questionnements d'identités développés dans le travail de C. Gilain sur le théorème fondamental de l'algèbre [Gilain 1991], celui de G. Cifoletti sur les pratiques algébriques de Pelletier et Gosselin dans le cadre de la « tradition algébrique française » de la renaissance [Cifoletti 1992] ainsi que l'étude des fluctuations d'élaborations mathématiques sur une longue période présentée dans l'ouvrage collectif coordonné par P. Benoit, K. Chemla et J. Ritter sur les fractions [Benoit et al. 1992] et le regard de C. Goldstein sur les relations entre analyse diophantienne et descente infinie [Goldstein 1993].

[3] Le terme « polynôme bilinéaire » d'abord utilisé par Jordan pour désigner l'expression $\Sigma A_{\alpha\beta} x_\alpha y_\beta$, est remplacé par l'expression « forme bilinéaire » dans le mémoire annoncé par la note de décembre 1873 et publié au printemps 1874. L'adoption du terme « forme » employé par Kronecker et qui renvoie à l'arithmétique de Gauss illustre la communication scientifique qui passe par la controverse. La question de la communication entre les deux centres que sont Paris et Berlin sera abordée de manière plus approfondie dans la suite de cet article.

[4] En des termes qui nous sont contemporains, la notion de forme bilinéaire joue pendant longtemps un rôle analogue à celui que jouera la notion de matrice dans l'algèbre linéaire du XX^e^ siècle. Dans les années 1874-1880, cette notion passe d'un sujet de recherche local à une théorie d'envergure internationale par sa grande étendue de domaines d'interventions comme la géométrie [Klein 1868], la théorie des formes quadratiques [Kronecker 1874] , [Darboux 1874], divers problèmes d'intégrations de systèmes différentiels linéaires comme les systèmes à coefficients constants [Jordan 1871 et 1872], les équations de Fuchs [Hamburger 1872] , [Jordan 1874c] ou le problème de Pfaff (Frobenius et Darboux 1875-1880).

[5] L'origine de la théorie des formes bilinéaires, à Berlin, dans les années 1860 est présentée dans [Hawkins 1977] et [Brechenmacher 2006a]. Au sujet des premiers travaux de Christoffel sur les formes consulter plus particulièrement [Mawhin 1981]. Pour des compléments mathématiques sur le problème de la réduction canonique des couples de matrices, voir [Dieudonné 1946].

l'examen du déterminant du « faisceau » (« Schaaren ») $A+sB$ et de l'équation polynomiale $|A+sB|=0$. Les racines $s_1, s_2, ..., s_n$ de l'équation caractéristique ne fournissaient à Kronecker un système complet d'invariants qu'à la condition que ces racines soient simples ([6]). La résolution générale du problème avait fait l'objet de la publication de deux mémoires en 1868. Le premier, de Weierstrass, introduisait un système d'invariants, les diviseurs élémentaires, obtenus par l'étude des systèmes successifs de sous déterminants du polynôme $|A+sB|$ et permettant la caractérisation des couples non singuliers de formes bilinéaires indépendamment de la multiplicité des racines caractéristiques ; le second, de Kronecker, traitait le cas singulier correspondant à l'occurrence d'un déterminant identiquement nul. Le théorème des diviseurs élémentaires s'imposait alors comme le résultat principal de la théorie des formes bilinéaires et Weierstrass montrait que la caractérisation des couples de formes donnait une même formulation à différentes questions abordées dans le passé comme l'intégration des systèmes d'équations différentielles linéaires ($AY''=BY$), la détermination des axes principaux des coniques et quadriques ($AY=\lambda Y$) ou la caractérisation des couples de formes quadratiques. Comme nous le verrons plus en détail dans la suite de cet article, c'est après avoir résolu à sont tour ces différents problèmes entre 1871 et 1872 à l'aide de méthodes élaborées dans le cadre de la théorie des groupes, que Jordan prit connaissance de la théorie des formes bilinéaires élaborée à Berlin et démontra que « la réduction simultanée de deux fonctions $P$ et $Q$ est un problème identique à celui de la réduction d'une substitution linéaire à sa forme canonique » [Jordan 1873, p. 1487].

Dans sa note de décembre 1873, Jordan intervient au cœur de la théorie des formes bilinéaires qu'il réorganise en développant une « méthode générale » consistant à « ramener » un « polynôme bilinéaire » à une « forme canonique simple » relative au type de substitutions considérées et sur laquelle il s'appuie pour distinguer trois types de questions :

> On sait qu'il existe une infinité de manières de ramener un polynôme bilinéaire
>
> $$P= \Sigma A_{\alpha\beta} x_\alpha y_\beta \ (\alpha=1,2,...,n, \beta=1,2,...,n)$$
>
> à la forme canonique $x_1y_1+...+x_my_m$, […] par des transformations linéaires opérées sur les deux systèmes de variables $x_1,..,x_n$, $y_1,..,y_n$. Parmi les diverses questions de ce genre que l'on peut se proposer, nous considérons les suivantes :
>
> 1. Ramener un polynôme bilinéaire $P$ à une forme canonique simple par des substitutions orthogonales opérées les unes sur $x_1,...,x_n$, les autres sur $y_1,...,y_n$.
>
> 2. Ramener $P$ à une forme canonique simple par des substitutions linéaires quelconques opérées simultanément sur les $x$ et les $y$.
>
> 3. Ramener simultanément à une forme canonique deux polynômes $P$ et $Q$ par des substitutions linéaires quelconques, opérées isolement sur chacune des deux séries de variables. [Jordan 1873, p. 1487].

La classification que nous pouvons voir à l'œuvre dans la citation ci-dessus est d'abord basée sur le type de groupe de substitutions agissant sur la ou les formes bilinéaires, elle présente la théorie des formes comme une mathématique pure dont l'organisation se rapproche de celle de la théorie des groupes réalisée par Jordan dans son traité de 1870. Dans le même temps, c'est sa capacité à s'appliquer dans divers domaines qui permet à la mathématique des formes de s'enrichir d'une épaisseur théorique ([7]). La « forme

[6] Par exemple la forme bilinéaire $B(X,U)=ux+uy+vy$ n'a qu'une seule valeur propre 1 qui ne suffit pas à la caractériser car $B$ n'est pas la forme identité.

[7] On dirait aujourd'hui que la « forme canonique » $x_1y_1+...+x_my_m$ permet de déterminer les classes d'équivalence des matrices carrées pour la relation d'équivalence ($ARB \Leftrightarrow \exists P,Q \in GL_n(\ ), PAQ=B$). Le problème 1. consiste en l'étude de la relation de similitude des matrices orthogonales ($ARB \Leftrightarrow \exists P \in O(\ ), P^{-1}AP=B$). Le problème 2 renvoie à la congruence des matrices ($ARB \Leftrightarrow \exists P \in GL_n(\ ), {}^tPAP=B$)). Le problème 3 à l'équivalence des couples de matrices $(A, B)$. Le problème 3 intervient pour la résolution des systèmes d'équations différentielles

canonique » $x_1y_1+...+x_my_m$ généralise la loi d'inertie de la théorie des formes quadratiques, le problème 1 fait référence à la classification des fonctions homogènes du second degré réalisée par Augustin Louis Cauchy en 1829 dans un cadre géométrique, le problème 2 renvoie à la question arithmétique de l'équivalence des formes quadratiques dans la tradition de Carl Friedrich Gauss et des *Disquitiones aritmeticae* de 1801 et le problème 3 provient de la théorie des systèmes d'équations différentielles linéaires à coefficients constants remontant au XVIII^e siècle. En présentant la théorie des formes bilinéaires comme donnant une identité nouvelle à des problèmes anciens issus de théories distinctes, la note de Jordan, qui est la première publication consacrée à cette théorie intervenant hors du centre berlinois, détache celle-ci du contexte local de sa création et la présente comme une théorie autonome dont il est licite de renouveler l'organisation et les méthodes. Ainsi, bien que la chronologie assure la priorité à ceux que Jordan désigne comme les « géomètres de Berlin », ce dernier revendique la nouveauté d'une méthode indissociable d'idéaux comme la simplicité et la généralité :

> Le premier de ces problèmes est nouveau, si nous ne nous trompons. Le deuxième a déjà été traité (dans le cas où *n* est pair) par M. Kronecker, et le troisième par M. Weierstrass ; mais les solutions données par les éminents géomètres de Berlin sont incomplètes, en ce qu'ils ont laissé de côté certains cas exceptionnels qui, pourtant, ne manquent pas d'intérêt. Leur analyse est en outre assez difficile à suivre, surtout celle de M. Weierstrass. Les méthodes nouvelles que nous proposons sont, au contraire, extrêmement simples et ne comportent aucune exception [...]. [*Ibidem*].

Dans sa réponse du 19 janvier 1874, Kronecker rejette aussi bien la nouveauté et la validité des méthodes de Jordan que l'organisation de la théorie autour de trois types de réductions canoniques. Il rappelle que, dès 1868, Weierstrass et lui-même avaient organisé la théorie des formes bilinéaires autour d'un unique problème, celui de la caractérisation des couples de formes.

> Dans le Mémoire de M. *Jordan [...],* la solution du premier problème n'est pas véritablement nouvelle ; la solution du deuxième est manquée, et celle du troisième n'est pas suffisamment établie. Ajoutons qu'en réalité ce troisième problème embrasse les deux autres comme cas particuliers, et que sa solution complète résulte du travail de M. *Weierstrass* de 1868 et se déduit aussi de mes additions à ce travail. Il y a donc, si je ne me trompe, de sérieux motifs pour contester à M. *Jordan* l'invention première de ses résultats, en tant qu'ils sont corrects [...]. [Kronecker 1874b, p. 1181 (les italiques sont dans le texte original)].

Durant l'hiver 1874, Kronecker détaille l'organisation qu'il donne à la théorie des formes bilinéaires lors de communications mensuelles à l'Académie de Berlin (annexe 1). De février à mars, ces interventions régulières restent cependant non polémiques vis-à-vis de Jordan. Au premier moment de la querelle de priorité que constituent les interventions publiques de décembre et janvier opposant deux organisations théoriques, deux méthodes et deux théorèmes succède, durant l'hiver, une correspondance privée entre les deux savants ([8]).

---

$AY''+BY=0$. Dans le cas particulier où $B=I$, la relation d'équivalence des couples $(A,I)$ est identique à la relation de similitude $B=P^{-1}AP$. Comme le fait remarquer Kronecker, le 3^e problème suffit à déduire les deux autres, le problème 1. revenant à l'étude de la congruence du couple $(A,I)$ et le 2. à l'équivalence du couple $(A, {}^tA)$.

[8] Cette correspondance, conservée dans les archives de l'école Polytechnique (côte VI2a2X1855), est éditée dans [Brechenmacher 2006a]. Voir plus généralement la présentation de la correspondance mathématique de Jordan par la conservatrice des archives, C. Billoux [Billoux 1985].

## 2. Un va-et-vient entre public et privé.

L'échange épistolaire débute lorsque Kronecker adresse à Jordan sa communication du 19 janvier [Kronecker 1874a] accompagnée d'une « sorte de sommation » de reconnaître publiquement la priorité des travaux berlinois sur la caractérisation des couples *(P, Q)* de formes bilinéaires ([9]) :

> ~~Notre ami commun~~ M. Kronecker qui vient ~~attaquer~~ de critiquer cette note avec une certaine vivacité dans les Monatsberichte, m'adresse en outre une sorte de sommation d'avoir à déclarer dans les Comptes Rendus 1° Que dans ses additions à votre travail de 1868 il a traité le cas $[P,Q]=0$ ; 2° que les formes canoniques que j'ai indiquées ne sont autres que la formule (44) de votre mémoire. (Jordan à Weierstrass, janvier 1874, dans [Brechenmacher 2006a, p. 50]).

La tentation d'une rhétorique guerrière dont témoignent les brouillons des lettres de Jordan atteste de la violence de la question de priorité :

> Je ~~ne voudrais pas~~ tiens en effet à certifier deux retards d'impression dont vous me faites part, ~~que je désire la guerre, et~~ que je préférerais ~~une polémique / guerre/~~ des débats publics à des explications amicales. Ce n'est pas moi qui ai ~~ouvert les hostilités~~ commencé la polémique. (Jordan à Kronecker, février 1874, dans [Brechenmacher 2006a, p. 64]).

La violence du conflit naît de sa dimension publique. Pour Jordan, c'est Kronecker qui « commence la polémique » par sa critique publique à l'Académie de Berlin. Pour Kronecker, Jordan est le premier fautif pour avoir publié sa note aux *Comptes rendus* sans, au préalable, être entré en contact avec lui-même ou Weierstrass :

> J'ai publié il est vrai (c'était mon droit évident) sans vous consulter des recherches qui complétaient les vôtres sur une question dont vous vous étiez occupé, et dont vous ne m'aviez jamais entretenu. ~~C'était mon droit évident.~~ Là-dessus, sans explication préalable à l'instant même vous publiez une critique plus longue que mon article, où vous me reprochez 1° De n'avoir rien compris à la manière de poser la question. 2° De n'y avoir apporté aucun élément nouveau 3° D'avoir pillé sans scrupule M. Weierstrass, M. Christoffel et vous. Si au lieu de jeter brusquement ce débat dans le public, vous vous étiez adressé à moi pour échanger des explications, comme je me voyais en droit de l'espérer, nous nous serions sans doute entendu. […] La publication imprévue de vos objections a un peu changé tout cela. ~~Il faut bien que la réponse soit publique~~ Attaqué devant tout le monde, il me faut bien répondre de même et il ne tiendra pas à moi que ce débat en reste là. (Jordan à Kronecker, janvier 1874, dans [Brechenmacher 2006a, p. 38]).

L'opposition public / privé est un des moteurs de la controverse et la correspondance entre les deux savants se présente comme une tentative de ramener le débat de la scène publique (journaux et communications académiques) à une relation privée. Dans la sphère privée, les témoignages de sympathie de Kronecker visent à apaiser le sentiment de Jordan d'avoir été agressé « devant tout le monde » : il ne s'agit ni d'une attaque personnelle (« je ne vous écrirais pas une lettre si longue et si précise »), ni d'une accusation de plagiat (« je pensais que mes "remarques" du 18 mai 1868 n'étaient pas connues de vous avant votre note aux comptes rendus » (Kronecker à Jordan, février 1874, dans [Brechenmacher 2006a, p. 54]). Quelle est exactement la dimension publique de la controverse ? Jordan redoute de faire l'objet d'une « réclamation collective » des « géomètres berlinois » dans un courrier adressé à Weierstrass dès janvier 1874 :

> Les termes ~~un peu ambigus~~ de la lettre de M. Kronecker me laissant incertain si je me trouve ici en présence d'une réclamation collective de votre part ou si je n'ai à répondre qu'à lui seul. J'aimerais être fixé à cet égard. [….] M. Kronecker a également relevé le

[9] Le terme de « sommation » employé par Jordan fait ici référence à un courrier de Kronecker de janvier 1874 qui n'a pas été retrouvé.

> passage où j'ai pris la liberté de parler des difficultés que j'ai éprouvées dans l'étude de votre mémoire. Je sais que mon impression à cet égard a été partagée par d'habiles géomètres que cette lecture a découragé. Je crois pouvoir attribuer une grande partie des difficultés à la forme synthétique que vous avez donné à votre démonstration. J'espère d'ailleurs que vous n'aurez pas été froissé de cette légère critique où il ne me semble pas avoir dépassé les bornes de la courtoisie. (Jordan à Weierstrass, janvier 1874, dans [Brechenmacher 2006a, p. 50]).

Cette lettre restera sans réponse et Weierstrass ne prendra pas part à la querelle. Des échanges avec Arthur Cayley, James Joseph Sylvester et Henry John Stephen Smith en prévision du voyage de Jordan à Londres au début de l'année 1875 attestent des tentatives du savant parisien d'obtenir des appuis publics. Ces tentatives restent vaines ([10]). La querelle semble d'ailleurs mettre Jordan en difficulté sur la scène parisienne - notamment vis-à-vis de Charles Hermite, le « cher ami » de Kronecker - et ruiner ses ambitions de remplacer Joseph Bertrand à l'Académie. Il n'y a pas d'engagement explicite des communautés savantes envers l'un ou l'autre des protagonistes. Si après la guerre de 1870, des tensions existent entre les géomètres parisiens et berlinois, comme en témoigne par exemple la correspondance de Joseph Liouville ([11]), pour ce qui est de la controverse entre Jordan et Kronecker, la scène publique est avant tout le lieu des débats mathématiques. La correspondance ne parvient en effet pas à une conciliation sur le plan mathématique, d'un côté Kronecker ne parvient pas, malgré les longs compléments de ses courriers, à justifier sa priorité de 1866 ; d'un autre, Jordan échoue dans sa tentative de convaincre de son « droit évident » à soutenir l'originalité de sa méthode comme ayant « plutôt mis en lumière que rabaissé l'importance » du résultat de Weierstrass « en montrant qu'il avait implicitement résolu l'un des problèmes fondamentaux de la théorie des substitutions linéaires, autrement féconde, à mon avis, que la théorie algébrique des formes du second degré » (Jordan à Kronecker, janvier 1874, dans [Brechenmacher 2006a, p. 38]). L'explication mathématique se tiendra, au printemps 1874, sur la scène publique des *Comptes rendus* de l'Académie de Paris.

## 3. La querelle : un révélateur de connaissances tacites et un vecteur de communication scientifique.

La correspondance révèle la difficulté pour Jordan d'acquérir la part tacite des pratiques développées à Berlin dans les années 1860 pour l'étude des formes bilinéaires. Jordan ne connaît ni les travaux de Christoffel, ni ceux de Kronecker de 1866 avant qu'un courrier ne le somme de les citer et, s'il révèle alors son ignorance, il en appelle à la vertu académique de ne pas revendiquer sur ce sujet :

> Je regrette surtout que vous ayez cru devoir introduire dans cette question le nom de M. Christoffel. […] J'ajouterai qu'à cette époque, je ne connaissais pas son mémoire ~~ce qui est bien permis.~~ Je suis persuadé d'ailleurs qu'il aurait été fort réservé dans les réclamations de ce genre. C'est d'après ce principe que j'ai laissé passer sans rien dire, dans le dernier numéro du Journal de M. Borchardt un mémoire volumineux de M. Sohnke sur la symétrie

[10] Jordan ne reçoit comme soutien direct qu'une remarque courte de Smith dans une lettre envoyée d'Oxford le 16 janvier 1875 : « I was wrong in saying that I could not follow <u>your</u> papers on the bilinear forms : it is M. Kronecker who gives me the trouble » [Archives Polytechnique VI2aX1855, lettre n°54].

[11] Une lettre de Kronecker est adressée à Liouville le 16 décembre 1875 à propos de la publication d'un mémoire sur la théorie arithmétique des formes quadratiques de déterminant négatif dans le journal de *Mathématiques pures et appliquées* (ayant été écarté de l'édition du journal, Liouville ne pourra pas publier ce mémoire). Dans sa lettre, Kronecker propose à Liouville et Chasles les deux places vacantes de correspondants étrangers de la section de géométrie de l'Académie berlinoise. La réponse de Liouville et ses annotations portées sur un brouillon de cours au collège de France témoignent de crispations consécutives à la guerre de 1870 entre les Académies de Paris et Berlin. Nous remercions Norbert Verdier pour avoir porté ces documents à notre attention.

> dans le plan où je n'étais pas cité, bien que j'eusse traité complètement cette question, pour le plan et pour l'espace, dans le journal de MM Brioschi et Cremona, il y a de cela cinq à six ans. (Jordan à Kronecker, janvier 1874, dans [Brechenmacher 2006a, p. 38]).

Il y a bien pour Jordan quelque chose d'étranger dans les pratiques associées aux formes bilinéaires. Il y a d'abord les difficultés éprouvées face à la lecture du mémoire de Weierstrass et confessées à de nombreuses occasions dans la correspondance. Au-delà des méthodes de calculs d'invariants basées sur la théorie du déterminant ([12]), les relations qu'entretiennent les publications berlinoises les unes avec les autres sont étrangères à Jordan. Ainsi, lorsqu'il publie sa note de 1873, celui-ci n'a pas conscience que le théorème de Weierstrass de 1868 vient répondre à des questions posées par Kronecker et Christoffel en 1866, et que si ce théorème se limite au cas des faisceaux non singuliers, il ne peut être dissocié du mémoire de Kronecker publié à sa suite dans le *Journal de Crelle*. Comme le précisera ce dernier, les deux publications de 1868 sont liées ensemble, comme un « développement en deux parties » (Kronecker à Jordan, mars 1874, dans [Brechenmacher 2006a, p. 68]). L'intitulé du mémoire de Kronecker ne faisant mention que des formes quadratiques, Jordan ne l'associe d'abord pas aux formes bilinéaires et se méprend sur son objet.

> Si au lieu de jeter brusquement ce débat dans le public, vous vous étiez adressé à moi pour échanger des explications, comme je me voyais en droit de l'espérer, nous nous serions sans doute entendu. Sur votre indication, j'aurais ~~constaté immédiatement, ce que j'ai reconnu trop tard, que votre méthode de 1868~~ relu plus attentivement votre mémoire de 1868 et constaté, ce que je n'avais pas remarqué à première vue, que les formes bilinéaires non citées dans votre travail, y sont pourtant implicitement comprises. […] Je me suis borné à dire en effet que vos solutions, dont je ne contestais nullement l'exactitude, laissaient de côté certains cas particuliers. Cela ne me paraît guère contestable pour la solution du second problème, contenue dans votre note du 15 octobre 1866, et quant au troisième problème, M. Weierstrass signale et exclut expressément dans son mémoire le cas où $[P,Q]=0$. […] Ne m'occupant que des formes bilinéaires, je n'avais pas à citer votre mémoire du 18 mai 1868, qui a trait aux formes quadratiques [….]. D'ailleurs vous n'aviez donné dans ce mémoire qu'un commencement de réduction et non une solution. […] Mais l'autre formule, relative au cas où $[P,Q]=0$, que vous passez sous silence, est précisément ce qu'il y a de nouveau dans ma note, ainsi que vous le reconnaitrez facilement en examinant la chose avec plus de soin. [*Ibidem*].

L'objet de ce mémoire, son articulation avec la publication de Weierstrass, sont des informations inaccessibles à Jordan et qui participent de la part tacite des pratiques d'une communauté, d'un réseau berlinois, dont la communication nécessite la longue lettre adressée par Kronecker à Jordan le 11 février 1874 :

> […] la formule finale que je donne aux faisceaux de formes quadratiques de déterminants identiquement nuls, associée au résultat de Weierstrass, donne une solution complète du problème de la transformation de deux faisceaux quelconques ([13]). Ils ne contiennent pas un « début » de solution mais la solution elle-même […] pour tous les faisceaux. (Kronecker à Jordan, février 1874, traduction F.B.(le terme souligné est dans le texte original), dans [Brechenmacher 2006a, p. 54]).

[12] Les méthodes de Weierstrass et Kronecker sont largement basées sur la théorie du déterminant. Kronecker en fait un des sujets principaux de ces cours d'université dans les années 1870-80. Selon le témoignage de Frobenius, Weierstrass utilise depuis 1864 une définition « quasi axiomatique » des déterminants dans laquelle [Knobloch, 1994] voit l'origine de la définition qui nous est contemporaine. La théorie du déterminant est essentielle dans les recherches sur les invariants de la seconde moitié du XIXe siècle, pour un exposé historique écrit au début du XXe siècle, consulter [Muir, 1906].

[13] Le terme allemand « Schaaren » employé par Kronecker pour désigner les expressions $uP+vQ$ ($P$, $Q$ deux formes bilinéaires ; $u$ et $v$ deux réels) est traduit ici par l'expression « faisceau » employée par Jordan.

Kronecker identifie très précisément la nature de la méprise de Jordan qu'il qualifie de « levis culpa », une faute légère et commune que Kronecker reconnaît avoir commis lui-même dans le passé vis-à-vis d'Hermite ou Joseph Alfred Serret. Réparer cette faute ne demande qu'une preuve de « moralité » et de « loyauté », que Jordan apportera en reconnaissant publiquement, au mois de mars 1874, sa méprise sur le contenu des travaux de Kronecker de 1868 ([14]). Que l'origine de la querelle soit liée à une difficulté de communication scientifique aurait pu permettre une sortie de crise rapide. Pourtant, la conciliation « sur cette même chose dont il n'y a plus rien à dire et qui me peine » échoue (Kronecker à Jordan, mars 1874, dans [Brechenmacher 2006a, p. 68]). Vecteur de communication scientifique, la correspondance a permis à Jordan de se familiariser avec les pratiques berlinoises relatives à l'étude des faisceaux de formes : méthodes de calculs d'invariants, chronologie des publications, architecture donnée à la théorie, interaction avec l'arithmétique des formes quadratiques. Les 2 mars et 27 avril 1874, les interventions publiques successives de Jordan et Kronecker à l'Académie de Paris, signalent l'échec de la « discussion privée », des « explications amicales » et la fin de la correspondance. A cette date, il ne s'agit plus de communiquer ou d'expliquer mais d'imposer une organisation de la théorie des formes bilinéaires.

## II. Une controverse opposant deux fins données à une histoire commune.

### 1. La généralité du théorème de Weierstrass perçue comme une rupture dans l'histoire des pratiques algébriques.

Comme nous le détaillons dans ce paragraphe, Kronecker associe au théorème de Weierstrass un idéal de généralité qu'il va expliciter par une critique des pratiques algébriques du passé et sur lequel il va s'appuyer pour dénoncer le caractère formel qu'il prête à la réduction canonique de Jordan.

> Ces problèmes, ainsi qu'ils sont posés ici [dans la note de Jordan de 1873], manquent tout à fait de précision même si le mot « canonique » suivant son sens propre donne justement l'impression qu'il pourrait s'agir de quelque chose absolument déterminé. […] La signification de l'expression « forme canonique » ou « forme canonique simple » utilisée par Mr. Jordan pour préciser la question [de la théorie des formes bilinéaires et quadratiques], n'a aucune pertinence générale ou décisive et désigne une notion sans aucun contenu objectif. [Kronecker 1874a, p. 367, traduction F.B.] ([15]).

Attribuant la même dénomination de forme canonique à trois expressions algébriques distinctes associées aux trois problèmes par lesquels il ordonne sa théorie des formes bilinéaires, Jordan est, dès janvier 1874, accusé par Kronecker de recourir à une notion sans signification précise, confondant « formel » et « contenu », « moyens » (formes canoniques) et « véritable objet de la recherche » (déterminer à quelles conditions une forme peut être transformée en une autre). Si Kronecker emploie lui-même des « formes normales » depuis 1866, et si ce que Jordan désigne comme « forme canonique » intervient dans la démonstration du théorème de Weierstrass de 1868 ([16]), de telles expressions algébriques ne sauraient quitter leurs « places relatives » de méthodes pour acquérir le statut d'une « notion » essentielle dans l'organisation donnée à la théorie. Plus généralement, la critique de Kronecker témoigne d'une conception selon laquelle le rôle de l'algèbre doit être du côté des méthodes, au service d'« autres disciplines » comme l'arithmétique qui, comme nous le verrons plus loin, doit pour le géomètre berlinois organiser la théorie des formes dans la tradition de Gauss :

> Si de telles expressions générales sont trouvées, on pourrait au besoin leur donner a posteriori une même désignation de formes canoniques pour des motivations de simplicité et de généralité. Mais si on ne veut pas en rester à ces aspects uniquement formels qui ont souvent été mis en avant par les travaux d'algèbre les plus récents -certainement pas au profit de la défense de la science–, alors on ne peut omettre de déduire le bien fondé de l'établissement des formes canoniques pour des raisons internes. En réalité, les dites formes canoniques ou formes normales sont effectivement déterminées uniquement par l'orientation donnée à l'étude et doivent donc seulement être considérées comme les moyens, mais non comme le but de la recherche. Cela ressort notamment clairement partout où le travail algébrique est effectué au service d'autres disciplines mathématiques, dont il reçoit ses fins et dont dépendent ses objectifs. Toutefois l'algèbre peut naturellement elle-même susciter également des motifs suffisants visant à l'établissement de formes canoniques, comme par exemple lorsque Mr. Weierstrass et moi-même avons été conduits dans les deux travaux cités par Mr. Jordan, à l'introduction de certaines formes normales, dont la place relative a été explicitement et clairement soulignée. [*Ibidem*, traduction F.B.] ([17]).

---

[14] Jordan fait observer que d'autres ont fait la même faute à son égard, voir la première citation de ce paragraphe.

[15] Sowie sie hier gestellt sind, ermangeln diese Probleme durchaus der Bestimmtheit, wie sehr auch grade das Wort « canonisch » seinem eigentlichen Sinne gemäss, den Schein von etwas absolut Bestimmtem zu erwecken geeignet ist. In der That hat der Ausdruck « canonische Form » oder « einfache canonische Form », welchen Hr. Jordan behufs Präcisirung der Frage gebraucht, keinerlei allgemein massgebende Bedeutung und bezeichnet an und für sich einen Begriff ohne jeden objectiven Inhalt.

[16] Consulter à ce sujet l'annexe 3. Sur la démonstration de Weierstrass de 1868 voir [Hawkins 1977] et [Brechenmacher 2006a].

La critique de Kronecker est précisée par la publication, en mai 1874, d'un appendice au mémoire de janvier. Le caractère formel de la notion de forme canonique est alors relié à la question de la « généralité » et aussi de l'« uniformité » de l'organisation donnée à la théorie des formes bilinéaires.

> Toutefois, il ne faut pas du tout être surpris que pour un développement à la fois tout à fait général et uniforme, comme [Jordan] en donne dans son travail cité, l'auteur soit nécessairement contraint de prouver certains nouveaux principes ; et il faudrait nous étonner au contraire, si conformément aux affirmations de Jordan (« Les méthodes nouvelles que nous proposons sont, au contraire extrêmement simples... » « On voit par une discussion très simple, que l'on peut transformer... ») les moyens les plus simples devraient suffire. [Kronecker 1874c, p. 404, traduction F.B.] ([18]).

Les revendications de simplicité de Jordan sont caricaturées comme un simplisme naïf qui en reste « à des aspects uniquement formels » pour parvenir à une présentation uniforme dont la généralité n'est qu'apparente et qui n'atteint pas la profondeur des « raisons internes » de la théorie ([19]). A partir d'un « défaut » relevé dans un calcul de Jordan, la mise au dénominateur d'une expression algébrique susceptible de s'annuler ([20]), Kronecker dénonce un certain type de pratiques algébriques qui s'appuient sur des expressions « prétendument générales » mais perdant en réalité toute signification dans certains « cas singuliers ». Le géomètre berlinois développe alors un discours sur deux significations opposées du terme généralité.

---

[17] Nachträglich, wenn dergleichen allgemeine Ausdrücke gefunden sind, dürfte die Bezeichnung derselben als canonische Formen allenfalls durch ihre Allgemeinheit und Einfachheit motiviert werden können ; aber wenn man nicht bei den bloss formalen Gesichtspunkten stehen bleiben will, welche –gewiss nicht zum Vortheil der wahren Erkenntnis- in der neueren Algebra vielfach in den Vordergrund getreten sind, so darf man nicht unterlassen, die Berechtigung der aufgestellten canonischen Formen aus inneren Gründen herzuleiten. In Wahrheit sind überhaupt die so genannten canonischen oder Normalformen lediglich durch die Tendenz der Untersuchung bestimmt und daher nur als Mittel, nicht aber als Zweck der Forschung anzusehen. Dies tritt namentlich überall da deutlich hervor, wo die algebraische Arbeit im Dienste andrer mathematischer Disciplinen geleistet wird und von ihnen Ausgangs – und Zielpunkt angewiesen erhält. Aber auch die Algebra selbst kann natürlich ausreichende Beweggründe zur Aufstellung canonischer Formen liefern, und so sind z.B. die Momente, welche Hrn. Weierstrass und mich in den beiden von Hrn. Jordan citirten Arbeiten bei Einführung gewisser Normalformen geleitet haben, an den bezüglichen Stellen klar und deutlich hervorgehoben.

[18] Dass sich aber für eine zugleich einheitliche und ganz allgemeine Entwickelung, wie sie in der ben erwähten Arbeit gegeben ist, gewisse neue Principien als nöthig erwiesen, kann durchaus nicht befremden, und es wäre im Gegentheil zu verwundern, wenn wirklich den *Jordan*'schen Behauptungen gemäss (« Les méthodes nouvelles que nous proposons sont, au contraire extrêmement simples... » « On voit par une discussion très simple, que l'on peut transformer... ») die allereinfachsten Mittel dazu ausreichen sollten.

> Car on est habitué –en particulier dans les questions algébriques – à trouver des difficultés largement nouvelles, si on veut se détacher de la restriction à ces cas, que l'on a coutume de désigner comme généraux.
>
> Aussitôt que l'on perce la surface de la prétendue généralité, excluant chaque particularité, on pénètre l'intérieur de la vraie généralité que toutes les singularités recouvrent, et l'on trouve généralement ainsi seulement les difficultés réelles de l'étude, ainsi que les abondants nouveaux points de vue et phénomènes qu'elle contient dans ses profondeurs. [*Ibidem*, traduction F.B.] ([21]).

Le théorème des diviseurs élémentaires de Weierstrass donne l'exemple parfait de la « vraie généralité » que Kronecker dépeint métaphoriquement comme un océan dont il faut percer la « surface » (la prétendue généralité) pour pénétrer les « profondeurs » en résolvant une « question algébrique dans ses moindres détails ». Par son traitement général et homogène de la caractérisation des couples de formes bilinéaires, il vient sanctionner les « résultats très insuffisants » de pratiques développées « durant tout un siècle » mais « de manière sporadique » et négligeant d'aborder les singularités en « n'osant pas faire tomber » la condition que le déterminant $S=|A+sB|$ contient des facteurs inégaux ([22]). Faisant implicitement référence aux travaux d'auteurs comme Joseph Louis Lagrange, Charles Sturm, Cauchy ou Carl Gustav Jacobi, ces pratiques, que nous préciserons dans le troisième paragraphe et auxquelles nous avons consacré une étude détaillée dans un article à paraître [Brechenmacher 200?], avaient été développées pour le cas symétrique des « systèmes quadratiques » dont elles exprimaient les solutions $x_i$ par des expressions polynomiales « générales » données en fonction de factorisations de $S$ et de ses sous déterminants $P_{1i}$ obtenus par développements par rapport à la $1^{re}$ ligne et $i^e$ colonne. L'expression générale (*) $\dfrac{\frac{P_{1i}}{S}(x)}{x-s_j}$ donnant $x_i^{s_j}=\dfrac{\frac{P_{1i}}{S}(s_j)}{x-s_j}$ où $(x_i^{s_j})_{1\le i\le n}$ désigne le système de solutions associé à la racine $s_j$ ([23]). Kronecker condamne l'utilisation qui avait été faite de telles

---

19 Les idéaux de simplicité de Jordan reposant sur des méthodes élaborées pour la théorie des groupes, la signification du terme « simplicité » est bien peu simple en elle-même. Nous verrons plus loin la relation entre l'idéal de simplicité de Jordan et une pratique algébrique de réduction issue des travaux des années 1860-1870 sur la résolubilité des équations algébriques.

20 Ce défaut, immédiatement corrigé par Jordan, est sans conséquence sur l'organisation théorique que celui-ci propose. Voir à ce sujet [Brechenmacher 2006a, p.689].

21 Denn man ist es gewohnt –zumal in algebraischen Fragen- wesentlich neue Schwierigkeiten anzutreffen, wenne man sich von der Beschränkung auf diejenigen Fälle losmachen will, welche man als die allgemeinen zu bezeichnen pflegt. Sobald man von der Oberfläche der sogenannten, jede Besonderheit ausschliessenden Allgemeinheit in das Innere der wahren Allgemeinheit eindringt, welche alle Singularitäten mit umfasst, findet man in der Regel erst die eigentlichen Schwierigkeiten der Untersuchung, zugleich aber auch die Fülle neuer Gesichtpunkte und Erscheinungen, welche sie in ihren Tiefen enthält.

22 C'est-à-dire la condition imposant des racines caractéristiques distinctes donc des matrices diagonalisables.

23 On dirait aujourd'hui que les coordonnées des vecteurs propres d'une matrice $A$ sont données par les colonnes non nulles de la matrice adjointe de la matrice caractéristique $A\text{-}xI$, c'est-à-dire la matrice des cofacteurs. Par exemple, pour la matrice :

$$A=\begin{pmatrix}1 & -1 & 0\\ -1 & 2 & 1\\ 0 & 1 & 1\end{pmatrix}$$

d'équation caractéristique et de mineurs : $S=-x(3-x)(1-x)$, $P_{11}(x)=(1-x)(2-x)-1$, $P_{12}(x)=(1-x)$ et $P_{13}=-1$.

expressions, susceptibles de prendre une valeur du type $\frac{0}{0}$ en cas d'occurrence de racines communes entre les équations obtenues par les sous déterminants successifs de $S=0$, avant que Weierstrass ne démontre en 1858, que toute racine d'ordre $p$ de $S$ est une racine d'ordre supérieur ou égal à $p-1$ pour $P_{1i}$, l'expression (*) étant par conséquent toujours parfaitement définie ([24]). Comme nous le verrons plus loin, le théorème de 1858 énonce que tout couple $(A = \sum_{i=1}^{n} A_{ij}x_i x_j \; B = \sum_{i=1}^{n} B_{ij}x_i x_j)$ de formes quadratiques - avec $A$ définie positive -, peut être écrit sous la forme $(\sum_{i=1}^{n} X_i^{\,2}, \sum_{i=1}^{n} s_i X_i^{\,2})$ où les racines $s_1, s_2, \ldots, s_n$ de $S=0$ sont toujours réelles, quelle que soit leurs multiplicités. Ce résultat généralise la loi d'inertie selon laquelle, comme l'énonçait Hermite, une forme quadratique homogène à $n$ variables $A$ peut s'écrire comme une somme de carrés de la forme $A=\Delta_{n-1}X_1^2+\frac{\Delta_{n-2}}{\Delta_{n-1}}X_2^2+\ldots+\frac{\Delta}{\Delta_1}X_n^2$, où $\Delta$, $\Delta_1, \Delta_2, \ldots, \Delta_{n-1}, 1$ sont les mineurs principaux successifs du déterminant de $A$. Comme le formule explicitement Darboux en 1874, cette loi se généralise aux couples $(A, I)$ envisagés comme des formes quadratiques polynomiales $A+sI$, les quotients de mineurs successifs qu'elle met en œuvre sont alors toujours définis car « si elle [l'équation caractéristique] a des racines multiples, une racine d'ordre $p$ devra annuler tous les mineurs d'ordre $p-1$ de l'invariant » [Darboux 1874, p. 367]. Le mémoire de Darboux est publié dans le *Journal de Mathématiques* à la suite du mémoire de [Jordan 1874a] qui est au cœur de la controverse avec Kronecker ; les deux mémoires se complètent et tous deux concernent des théorèmes de Weierstrass, d'une part Jordan propose une nouvelle théorie des faisceaux de formes

---

L'expression $((1-x)(2-x)-1, 1-x, -1)$ donne une écriture polynomiale des coordonnées d'un vecteur propre dont $|S|$ donne le carré de la norme. Pour la valeur propre $s_1=1$, $\frac{S}{x-1}=x(3-x)$. Par conséquent,

$$x_1^{s_1} = \frac{P_{11}}{\frac{S}{x-1}}(1) = -\frac{1}{2}, x_2^{s_1} = \frac{P_{12}}{\frac{S}{x-1}}(1) = 0, x_3^{s_1} = \frac{P_{11}}{\frac{S}{x-1}}(1) = -\frac{1}{2},$$

On en déduit les coordonnées d'un vecteur propre normé associé à la valeur propre 1 : $(1/\sqrt{2}, 0, 1/\sqrt{2})$.
En procédant de manière analogue pour les valeurs propres $s_2=0$,

$$x_1^{s_2} = \frac{1}{3}, x_2^{s_2} = \frac{1}{3}, x_3^{s_2} = -\frac{1}{3}$$

et $s_3=3$,

$$x_1^{s_3} = \frac{1}{6}, x_2^{s_3} = \frac{2}{6}, x_3^{s_3} = -\frac{1}{6}$$

on obtient respectivement pour vecteurs propres normés, $(1/\sqrt{3}, 1/\sqrt{3}, -1/\sqrt{3})$ et $(1/\sqrt{6}, 2/\sqrt{6}, -1/\sqrt{6})$.
Interprétation dans le cadre des formes quadratiques : la forme associée à $A$ dans la base canonique de À$^3$,

$$A(x_1, x_2, x_3) = x_1^2 - 2x_1x_2 + x_2^2 + 2x_2x_3 + x_3^2 = 1.X_1^2 + 0.X_2^2 + 3.X_3^2$$

avec $(X_1, X_2, X_3) = (x_i^{s_j})_{i,j=1,2,3} A (x_i^{s_j})^{-1}_{i,j=1,2,3}$ où $(x_i^{s_j})_{i,j=1,2,3}$ est la matrice orthogonale de passage (ou de changement de base orthonormées) :

$$(x_i^{s_j})_{i,j=1,2,3} = \begin{pmatrix} \frac{1}{\sqrt{2}} & \frac{1}{\sqrt{3}} & \frac{1}{\sqrt{6}} \\ 0 & \frac{1}{\sqrt{3}} & \frac{2}{\sqrt{6}} \\ \frac{1}{\sqrt{2}} & -\frac{1}{\sqrt{3}} & -\frac{1}{\sqrt{6}} \end{pmatrix}$$

Interprétation dans le cadre des formes bilinéaires (symétriques) :

$$A(x,y) = x_1y_1-x_1y_2-y_1x_2+x_2y_2+x_2y_3+y_3x_2+x_3y_3=1.X_1Y_1+0.X_2Y_2+3.X_3Y_3.$$

[24] Contrairement au cas bilinéaire, dans le cas quadratique étudié par Weierstrass en 1858 les matrices sont symétriques donc toujours diagonalisables sur À. Une matrice est diagonalisable si et seulement si ses diviseurs élémentaires sont simples.

bilinéaires abordés par Weierstrass en 1868, d'autre part Darboux propose une nouvelle démonstration du résultat de 1858 sur les faisceaux de formes quadratiques ([25]).

> Ceci se confirme partout dans les rares questions algébriques qui sont mises en œuvre complètement jusqu'à leurs moindres détails, notamment dans la théorie des faisceaux des formes quadratiques qui a été développée plus haut dans ses caractéristiques principales. Parce que, pendant si longtemps, on n'osait pas faire tomber la condition que le déterminant ne contient que des facteurs inégaux, on est arrivé avec cette question connue de la transformation simultanée de deux formes quadratiques; qui a été si souvent traitée depuis un siècle, mais de manière sporadique, à des résultats très insuffisants et les vrais aspects de l'étude ont été ignorés. Avec l'abandon de cette condition, le travail de Weierstrass de l'année 1858 a conduit à un aperçu plus élevé et notamment à un règlement complet du cas, dans lequel n'existent que des diviseurs élémentaires simples. Mais l'introduction générale de cette notion de diviseur élémentaire, dont seule une étape provisoire était alors accomplie, intervient seulement dans le mémoire de Weierstrass de l'année 1868, et une lumière tout à fait nouvelle est ainsi faite sur la théorie des faisceaux pour n'importe quel cas, avec la seule condition que le déterminant soit différent de zéro. Quand j'ai aussi dépouillé cette dernière restriction et l'ai développé à partir de la notion de diviseur élémentaire des faisceaux élémentaires généraux, la clarté la plus pleine s'est répandue sur une quantité de nouvelles formes algébriques, et par ce traitement complet de l'objet des vues plus élevées ont été acquises sur une théorie des invariants comprise dans sa vraie généralité. [*Ibid.*, traduction F.B.] ([26]).

La querelle mathématique voit Kronecker porter un regard historique sur des pratiques algébriques du passé que l'on peut qualifier de « génériques » pour reprendre le qualificatif attribué par Thomas Hawkins, à l'occasion de ses travaux consacrés à l'origine de la théorie spectrale des matrices, à un certain type de raisonnements qui ne se préoccupe pas de la signification des symboles [Hawkins 1977, p. 122]. Déjà en germe dans le résultat de 1858 sur les formes quadratiques qui montre que la question essentielle n'est pas la multiplicité des racines de $S=0$ mais les relations de divisibilité entre le déterminant caractéristique et ses mineurs successifs (annexe 3), le théorème des diviseurs élémentaires est considéré comme rompant avec de telles pratiques. En effet, ce théorème a un caractère à la fois rigoureux et homogène rigoureux, par opposition aux formules génériques qui perdent toute signification en cas d'occurrence de racines multiples, homogène par opposition à une perception des différents cas de multiplicité de racines comme autant de cas singuliers nécessitant des traitements spécifiques. Cette rupture se manifeste notamment dans l'évolution des travaux de Kronecker. Dans un mémoire consacré aux fonctions thêta en 1866, celui-ci avait été amené à considérer l'occurrence de racines caractéristiques multiples comme un cas singulier nécessitant un traitement spécifique distinct de la méthode algébrique développée pour ce qu'il désignait alors comme le « cas général ». La résolution homogène permise par le système complet d'invariant introduit par Weierstrass en 1868 avait démontré la possibilité de concilier une approche algébrique et un traitement général du problème. Comme l'a montré Thomas Hawkins, elle était à l'origine de nouveaux idéaux de généralité dans les raisonnements algébriques se manifestant notamment dans l'élaboration de la théorie des formes bilinéaires à Berlin [Hawkins 1977, p. 157]. Remarquer que les premières lignes de l'histoire qui sanctionne des pratiques du passé comme « génériques » sont écrites par Kronecker lui-même dans le contexte très particulier d'une controverse mathématique ([27]), nous amène à suivre le regard différent porté par Jordan sur ces mêmes pratiques ainsi que la question de leur généralité et, par là, à aborder ce que le savant berlinois désigne comme « l'histoire de la théorie des faisceaux de formes quadratiques depuis un siècle » selon une approche complémentaire à celle suivie par Hawkins lorsqu'il avait étudié l'évolution du caractère générique des raisonnement de l'algèbre, de leur origine dans l'analyse de François Viète au développement d'exigences de rigueur au XIX^e siècle.

---

[25] Le mémoire de Darboux concerne avant tout les surfaces, et sa méthode est insérée par Gundelfinger dans la troisième édition de la géométrie analytique de Hesse. Pour une description détaillée du travail de Darboux, voir l'article de Meyer et Drach dans *L'encyclopédie des sciences mathématiques* de 1907.

[26] Diess bewährt sich durchweg in den wenigen algebraischen Fragen, welche bis in alle ihre Einzelheiten vollständig durchgeführt sind, namentlich aber in der Theorie der Schaaren von quadratischen Formen, die obein in ihren Hauptzügen entwickelt worden ist. Denn so lange man es nicht wagte, die Voraussetzung fallen zu lassen, dass die Determinante nur ungleiche Factoren enthalte, gelangte man bei jener bekannten Frage der gleichzeitigen Transformation von zwei quadratischen Formen, welche seit einem Jahrhundert so vielfach, wenn auch meist blos gelegentlich, behandelt worden ist, nur zu höchst dürftigen Resultaten, und die wahren Gesichtpunkte der Untersuchung blieben gänzlich unerkannt. Mit dem Aufgeben jener Voraussetzung führte die *Weierstrass*'sche Arbeit vom Jahre 1858 schon zu einer höheren Einsicht und namentlich zu einer vollständigen Erledigung des Falles, in welchem nur einfache Elementartheiler vorhanden sind. Aber die allgemeine Einführung dieses Begriffes der Elementartheiler, zu welcher dort nur ein vorläufiger Schritt gethan war, erfolgte erst in der *Weierstrass*'schen Abhandlung vom Jahre 1868, und es kam damit ganz neues Licht in die Theorie der Schaaren für den Fall beliebiger, doch von Null verschiedener Determinanten. Als ich darauf auch diese letzte Beschränkung abstreifte und aus jenem Begriffe der Elementartheiler den allgemeineren der elementaren Schaaren entwickelte, verbreitete sich die vollste Klarheit über die Fülle der neu auftretenden algebraischen Gebilde, und bei dieser vollständigen Behandlung des Gegenstandes wurden zugleich die wertvollsten Einblicke in die Theorie der höheren, in ihrer wahren Allgemeinheit aufzufassenden Invarianten gewonnen.

## 2. Une incorrection dans une pratique remontant à Lagrange.

Dans une note adressée à l'Académie de Paris en 1870, l'astronome Antoine Yvon-Villarceau appelle « l'attention des géomètres » sur un « point assez important de la théorie des équations linéaires », une « incorrection » dans la méthode « d'intégration des équations différentielles du mouvement de rotation d'un corps solide, soumis à l'action de la pesanteur [...], présentée pour la première fois par l'illustre auteur de la Mécanique Analytique, dans le cas des petites oscillations » d'une corde fixée en un point, lestée d'un nombre quelconque de masses et écartée de sa position d'équilibre [Yvon-Villarceau 1870, p. 762]. Le principe de conservation des forces vives permet de mathématiser le problème par un système d'équations différentielles linéaires à coefficients constants ([28]).

$$(a)\begin{cases} g\dfrac{d^2u}{dt^2}+a\dfrac{d^2s}{dt^2}+cu=0, \\ f\dfrac{d^2s}{dt^2}+a\dfrac{d^2u}{dt^2}+cs=0, \end{cases}$$

La méthode d'intégration du système *(a)* repose sur la détermination, par les méthodes d'éliminations, d'une équation algébrique

$$\frac{c^2}{\rho^4}-(f+g)\frac{c}{\rho^2}+fg-a^2=0$$

Aux racines $\rho$ et $\rho'$ de cette équation de degré 2 sont associées deux équations différentielles indépendantes, $\dfrac{d^2u}{dt^2}+\rho u=0$ et $\dfrac{d^2u}{dt^2}+\rho' u=0$, auxquelles est ramenée l'intégration du système.

---

[27] Cette perception de l'histoire s'avère donc indissociable de l'énoncé du théorème des diviseurs élémentaires que Kronecker oppose au théorème de réduction canonique de Jordan. Sur la question de la fabrication de l'histoire par les textes mathématiques, voir la discussion de référence de [Goldstein 1995], celle de [Dhombres 1998] en relation avec le concept de postérité, les exemples donnés par [Cifoletti 1992 et 1995] et, pour notre propos, les développements proposés dans [Brechenmacher 2006c et 200?].

[28] Dans la citation de Villarceau, *u* et *s* sont des fonctions de *t*, *g*, *f* et *a* des constantes. L'intervention de *u* et *s* est en miroir dans les équations, le système est donc symétrique. Dans une publication à paraître nous montrons que cette propriété se dégage d'une pratique algébrique spécifique élaborée par Lagrange en 1766 [Brechenmacher 200?].

$$(f)\begin{cases} s = \alpha\sin(\rho t+\beta)+\alpha'\sin(\rho' t+\beta') \\ u = \dfrac{a\rho^2}{c-g\rho^2}\alpha\sin(\rho t+\beta)+\dfrac{a\rho'^2}{c-g\rho'^2}\alpha'\sin(\rho' t+\beta') \end{cases}$$

La réduction du système à deux équations indépendantes nécessite donc l'obtention de deux racines distinctes et Yvon-Villarceau critique l'interprétation donnée par Lagrange selon laquelle l'inégalité des racines serait garantie par l'hypothèse de la stabilité des petites oscillations. En cas de racines multiples, les oscillations ne resteraient pas bornées car le temps *t* « sortirait du sinus » et les solutions prendraient la forme $s=t\sin(\rho t+\beta)$([29]).

> Je dis qu'il n'est pas nécessaire que cette condition soit remplie, pour que les petites oscillations se maintiennent. [...] Voici un cas très simple, auquel correspondent des racines égales de l'équation caractéristique : c'est celui d'un corps solide, homogène et de révolution, oscillant autour d'un point pris sur son axe de figure. Chacun comprendra sans recourir au calcul, que la petitesse des oscillations est assurée dans ce cas, si le centre de gravité est, à l'origine du mouvement, au-dessous du centre de suspension, à une petite distance de la verticale passant par ce point, et si le mouvement oscillatoire initial est suffisamment faible.[Yvon-Villarceau 1870, p. 765].

Si, comme Kronecker, Yvon-Villarceau met en évidence un « défaut » dans une méthode centenaire, défaut relatif à la négligence d'un cas singulier, le discours ne vise pas le caractère générique du raisonnement mais la représentation mécanique sous jacente à celui-ci : un mouvement oscillatoire doit pouvoir se décomposer en des mouvements indépendants. Cette intervention, qui a pour objet de mettre en évidence la possibilité de concilier des solutions stables, une décomposition en mouvements indépendants et l'existence de périodes propres égales, est issue de préoccupations mécaniques concernant notamment l'application de la méthode des petites oscillations aux mouvements séculaires des planètes sur leurs orbites ([30]). Remettant en cause une pratique consistant à donner une représentation mécanique (les solutions ne sont pas stables) à une propriété algébrique (la multiplicité des racines) d'une équation particulière (l'équation caractéristique) ([31]), la critique de l'astronome concerne les « géomètres » et la « théorie des équations linéaires » car elle pose la question théorique de la caractérisation des systèmes se ramenant à des équations indépendantes.

> J'ai cru devoir appeler l'attention des géomètres sur un point assez important de la théorie des équations linéaires, et qui n'occupe pas une place suffisante dans les traités sur cette matière. Peut-être la question que je soulève a-t-elle déjà été résolue; mais il faut croire que la solution n'est pas généralement connue, puisque l'incorrection que je signale dans la *Mécanique analytique* a pu échapper à un géomètre aussi érudit que le savant auteur de la nouvelle édition d'un ouvrage devenu classique ([32]). [*Ibidem*, p.766].

Jordan répond à l'appel de Villarceau par la publication de deux notes aux *Comptes rendus*. Il donne d'abord, en 1871, une méthode permettant l'intégration des systèmes différentiels à coefficients constants indépendamment de la multiplicité des racines : le système est associé à une substitution linéaire dont la réduction canonique permet de donner au système une forme intégrable [Jordan 1871, p.787]. Jordan donne aussi une caractérisation des systèmes se « réduisant » à des « équations distinctes » de la forme $\frac{dy_1}{dt} = \sigma y_1$ par une relation entre le déterminant caractéristique et ses mineurs : chaque racine caractéristique $\sigma$ de multiplicité $\mu$ doit être une racine des mineurs d'ordre $\mu$-1 extraits du déterminant caractéristique ([33]). Dans une seconde note publiée en 1872 et intitulée « Sur les oscillations infiniment petites des systèmes matériels », il répond plus précisément à la question posée par Villarceau en abordant le cas des systèmes intervenant en mécanique et démontre que le caractère quadratique de ces systèmes implique la relation entre déterminants et mineurs qui garantie la possibilité de les réduire à des « équations distinctes ». Ce faisant, Jordan dépasse le seul problème de l'intégration des systèmes symétriques en se référant au cadre plus général de la théorie des formes quadratiques et aux travaux de Cauchy et Hermite :

> Il est clair que la question de la réduction du système (1) à la forme canonique (7) est identique à ce problème connu : *Faire disparaître les angles des variables à la fois dans les deux formes quadratiques* T *et* U. [Jordan 1872, p. 320].

[29] Le système est stable et se ramène à des « équations distinctes qui s'intègrent isolément » s'il est diagonalisable dans À. Une matrice dont les valeurs propres sont toutes distinctes est diagonalisable mais la condition n'est pas nécessaire, une matrice présentant des valeurs propres multiples peut être diagonalisable ou non diagonalisable. La matrice du système d'équations différentielles linéaires à coefficients constants de Lagrange est toujours diagonalisable car symétrique.

[30] Les travaux de mécanique céleste d'Yvon-Villarceau concernent la détermination des orbites de divers corps célestes. Des exemples célèbres sont le calcul de la périodicité de la comète d'Arrest (1851) et la prévision des éphémérides en tenant compte des perturbations produites par Jupiter. Consulter à ce sujet [Baillaud 1957].

[31] Pour Yvon-Villarceau, la forme $s=t\sin(\rho t+\beta)$ serait attribuée aux solutions en cas de racine double par une généralisation abusive de l'intégration d'une équation linéaire d'ordre *n* à celle d'un système de *n* équations.

[32] Il s'agit ici d'une allusion à l'édition de la *Mécanique Analytique* par Bertrand en 1853.

Les propriétés par lesquelles Jordan caractérise, en réponse à l'appel de Villarceau, les systèmes réductibles à des équations distinctes sont en fait identiques à celles déjà énoncés par Weierstrass en 1858 (cas symétrique des couples de formes quadratiques) et 1868 (cas général des couples de formes bilinéaires) ([34]). Cette identité est d'abord relevée par Meyer Hamburger qui, en 1873, mêle les méthodes de réduction canonique de Jordan aux résultats de Weierstrass pour résoudre le problème posé par l'occurrence de racines caractéristiques multiples dans l'intégration des équations différentielles linéaires à coefficients non constants, dites équations de Fuchs [Hamburger 1873, p. 113]. La parution du mémoire d'Hamburger attire l'attention de Jordan sur la théorie des formes bilinéaires développée à Berlin et sur laquelle ce dernier va publier une série de notes entre 1873 et 1874. Jordan démontre que la transformation des formes peut s'interpréter comme l'action de groupes de substitutions (linéaires, orthogonales etc.) et, comme nous l'avons vu, prouve que le théorème des diviseurs élémentaires des couples de formes bilinéaires $(P,Q)$ « se ramène identiquement à celui de la réduction des substitutions linéaires à leur forme canonique ».

> Soient $\omega$ et $\omega'$ deux constantes quelconques, telles que les polynômes
>
> $$\mathcal{P} = \omega P+Q, \quad \mathcal{L}=\omega' P+Q$$
>
> Aient leurs déterminants différents de zéro. On exprimera aisément $P$ et $Q$ en fonction de $\mathcal{P}$ et $\mathcal{L}$. Reste à assigner une forme simple à ces deux derniers polynômes. [Jordan 1874a, p. 52].

Choisissant deux valeurs $\omega$ et $\omega'$ de la variable et posant $\mathcal{P}=\omega P+Q$ et $\mathcal{L}=\omega' P+Q$, Jordan réduit la forme bilinéaire $\mathcal{P}$ à sa « forme canonique » et désigne par $S$ la substitution permettant de transformer la forme bilinéaire $\mathcal{P}$ en $\mathcal{L}$ ([35]). Les propriétés de cette substitution $S$ « correspondante » au couple de formes $(P, Q)$ permettent d'expliciter l'identité mathématique des problèmes de transformations de formes et de réduction des substitutions : l'opération sur « les deux séries de variables » qui laisse la forme $\mathcal{P}$

[33] Cette propriété caractérise les endomorphismes diagonalisables par la condition que le degré de multiplicité $\mu$ d'une valeur propre doit coïncider avec la dimension du sous espace propre associé.

[34] Il semble que les travaux de Weierstrass de 1858 aient influencé Jordan de manière indirecte, par l'intermédiaire d'un mémoire publié à l'Académie de Saint Petersbourg par J. Somoff, porté à l'attention de Jordan par Yvon-Villarceau en 1872 et contenant un exposé de la méthode de Weierstrass appliquée au problème des petites oscillations (J. Somoff, Sur l'équation algébrique à l'aide de laquelle on détermine les oscillations très-petites d'un système de points matériels, *Mém. Acad. Sci. St. Pétersbourg,* (7) 1 no.14 (1859)).

[35] En des termes qui nous sont contemporains, si

$$\mathcal{P}={}^tXY, \text{ alors } \mathcal{P}.S=\mathcal{L}={}^tXSY,$$

Cela revient à faire opérer à droite le groupe linéaire sur l'ensemble des formes bilinéaires

$$\{\omega P+Q / \omega \in À \ , |\omega P+Q| \neq 0\}.$$

invariante correspond à l'opération de conjugaison des substitutions linéaires, $S \rightarrow U^{1}SU$, définie dans le *Traité des substitutions* de 1870 ([36]). La question de la transformation des couples de formes « se réduit » alors à celle de la composition des substitutions et il suffit, pour obtenir la réduction simultanée des formes $\mathcal{P}$ et $\mathcal{L}$ à leurs formes canoniques, de choisir la substitution $U$ de telle sorte que la substitution $USU^{-1}$ associée à $\mathcal{L}$ soit réduite à sa forme canonique ([37]).

### 3. La discussion sur l'équation à l'aide de laquelle on détermine les inégalités séculaires des planètes (1766-1874).

Comme nous l'avons vu, une première identité entre les théorèmes de Jordan et Weierstrass s'était manifestée, entre 1870 et 1874, par la capacité de ces deux énoncés à donner une résolution générale à de mêmes problèmes abordés dans le passé par des auteurs comme Lagrange, Cauchy ou Hermite dont les résultats étaient apparus lacunaires car limités à ce qui était désormais considéré comme un cas singulier, l'occurrence de racines caractéristiques distinctes. En ce sens, la controverse de 1874 oppose deux théorèmes se présentant comme deux clôtures différentes d'une histoire commune. Afin d'éclairer le rôle joué par cette histoire dans la querelle, nous avons mis en œuvre une recherche basée sur les références données par Jordan et Kronecker et complétée par un épuisement systématique des indications bibliographiques jusqu'à obtenir un corpus regroupant des textes mathématiques en un ensemble cohérent et limité dans le temps (1766-1874) ([38]). Quelques nœuds apparaissent dans l'enchevêtrement des références bibliographiques visualisé par le graphe simplifié de l'annexe 4. En 1766, Lagrange travaille sur l'intégration des systèmes d'équations différentielles linéaires à coefficients constants ; en 1829, la question traitée par Cauchy est la classification des surfaces du second degré et, en 1858, il s'agit pour Weierstrass de caractériser les transformations de fonctions homogènes du second degré par substitutions sur les variables. Notre corpus ne se laisse pas aisément caractériser. Ce n'est pas une théorie si on ne le fixe pas à un moment précis et qu'on le considère sur la longue durée. Ce n'est pas non plus le développement de la résolution d'un problème, nous verrons que s'il y a bien un problème fondateur, celui-ci est considéré, jusqu'au mémoire de Weierstrass de 1858, comme résolu par Lagrange dès 1766. En tenant compte des citations et références que les textes entretiennent entre eux, le corpus articule une discussion entre différents auteurs et différentes théories que nous désignerons sous le nom de discussion sur l'équation à l'aide de laquelle on détermine les inégalités séculaires des planètes ou, plus rapidement, discussion sur l'équation des petites oscillations. En envisageant notre corpus comme une discussion, nous voulons questionner son identité et cette question nous amène à nous démarquer de l'historiographie qui s'est intéressée aux textes de notre réseau dans l'objectif de dresser l'histoire d'une théorie, la théorie spectrale des matrices, permettant de formuler les différents problèmes étudiés comme une même question de réduction d'un couple de matrices *(A,B)* en un couple (*D, I*) où *A* est une matrice symétrique, *B* une matrice symétrique définie, *D* une matrice diagonale et *I* la matrice identité ([39]). La discussion jouant, entre 1870 et 1874, le rôle d'une histoire commune au travers de laquelle de premières identités se manifestent entre les travaux de Jordan sur les substitutions et de Kronecker et Weierstrass sur les formes bilinéaires, elle semble bénéficier d'une cohérence propre préalable à l'identité théorique donnée rétrospectivement par des théories algébriques comme la théorie des matrices des années 1930 ou la théorie des formes bilinéaires de 1874. Nous allons préciser cette identité en examinant les principaux nœuds de notre réseau que sont les travaux de Lagrange, Cauchy et Weierstrass.

---

[36] Si $T$ est une substitution « opérée sur les variables $x_1,\ldots,x_n$ » (opération à gauche) alors la substitution correspondante sera ${}^tTS$ car $T.\mathcal{L}={}^t(LX)SY={}^tX{}^tTSY={}^tTS.\mathcal{P}$ ). Si $T$ est une substitution opérant sur les « variables » $x_i$ (à gauche) et si $U$ est une substitution opérant sur les variables $y_i$ (à droite) alors la substitution correspondante sera ${}^tTSU$. Pour que l'opération de $T$ et $U$ n'altère pas la forme de $\mathcal{P}$, il faut que ces substitutions soient inverses l'une de l'autre, ${}^tT=U^{-1}$.

[37] La forme $\mathcal{P}$ est réduite à l'identité et la forme $\mathcal{L}$ à une forme dite de Jordan des formes bilinéaires. Jordan ne donne cependant pas explicitement la forme réduite des couples de formes et se contente de ramener le problème à la réduction des substitutions.

[38] Notre méthodologie est donc basée sur le choix de la querelle de 1874 comme moment de référence. D'autres choix de dates et de contextes feraient apparaître des réseaux différents ou jouer des rôles importants à des textes apparaissant ici comme des références mineures. Un autre point de vue, très proche puisque provenant des années 1880 et d'auteurs comme H. Poincaré, conduirait à ouvrir le corpus à des travaux d'astronomie et à centrer le questionnement sur la stabilité du système du monde et le problème des trois corps (voir à ce sujet la récente thèse de doctorat d'Anne Robadey [Robadey 2006) et, pour une perspective sur le long terme, [Laskar 1992]). En prenant pour moment de référence la date de 1890 et les travaux d'E. Weyr sur les matrices, nous serions amené à donner une plus grande importance aux travaux dans lesquels J. J. Sylvester définit, en 1851, le terme « matrice » mais qui ne jouent ici qu'un rôle mineur et ne sont cités qu'une seule fois par G. Darboux en 1874 (voir [Brechenmacher 2006d]). Un point de vue encore différent, des années trente du XX$^{e}$ siècle, pourrait faire participer d'un même corpus des textes de Lagrange et de C. Fourier, et y voir une origine de l'algèbre linéaire (voir les commentaires historiques des traités de [Mac Duffee 1933] et [Aitken et Turnbull 1932]).

Comme l'illustre l'intervention d'Yvon-Villarceau que nous avons déjà évoqué, la référence systématique à la résolution par Lagrange du « problème des oscillations très petites d'un système quelconque de corps » est une première caractéristique du corpus. Il s'agit d'étudier le déplacement d'un système de corps *m, m', m'', &c,* de coordonnées $\xi, \psi, \varphi, ...$, dont on suppose qu'il s'éloigne très peu d'une position d'équilibre. Les équations du mouvement se présentent sous « une forme linéaire avec des coefficients constants ». L'intégration du système est basée sur le principe de réduction de l'ordre d'une équation différentielle par la connaissance de solutions particulières, trouver *n* solutions particulières $\xi=E\sin(t\sqrt{K}+\varepsilon)$ de *n* équations indépendantes $\frac{d^2\xi}{dt^2}+K\ \xi=0$ permet, « en les joignant ensemble », d'exprimer les solutions générales, les valeurs de *K* apparaissant, par usage des méthodes d'élimination, comme les racines d'une équation résultante de degré *n*. La question se ramène donc à la résolution d'une équation algébrique qui s'avère caractéristique du problème des petites oscillations. Pour Lagrange, la nature mécanique des oscillations, supposées petites dans le problème posé, implique en effet un caractère spécifique de la nature algébrique des racines de l'équation, celles-ci doivent être réelles, négatives et inégales ([40]) :

[39] Voir par exemple [Gantmacher 1959, p. 311] pour des compléments mathématiques. Une présentation détaillée de la discussion des petites oscillations sur toute sa durée fera l'objet d'une publication ultérieure [Brechenmacher 200?] proposant une approche complémentaire aux travaux menés par l'historien T. Hawkins de l'université de Boston qui, entre 1974 et 1977, a consacré plusieurs articles à l'histoire de la théorie spectrale des matrices et en a fixé la structure qui a ensuite été reprise par les traités d'histoire des mathématiques proposant une histoire de l'algèbre linéaire. Une première étape correspond à une « origine » de la théorie spectrale dans les travaux de mécaniques du XVIII[e] siècle impliquant des systèmes différentiels linéaires à coefficients constants du type $AY'' = BY$ et dont l'intégration repose sur la recherche de valeurs propres $\lambda_i$ telles que $AX=\lambda_i BX$ permettant d'écrire le système sous la forme $DY'' = Y$ où $D=(\lambda_1, \lambda_2, ..., \lambda_n)$ est une matrice diagonale. A cette première étape succède un développement théorique initié par [Cauchy 1829] et complété par [Weierstrass 1858] se caractérisant par la démonstration de la nature réelle des valeurs propres des matrices symétriques. La troisième étape correspond à l'organisation de la théorie des formes bilinéaires, fixée par [Frobenius 1878] autour du théorème des diviseurs élémentaires énoncé dans [Weierstrass 1868] et permettant de caractériser les classes d'équivalence des polynômes de matrices (et donc les classes de similitudes).

[40] Dans cette conclusion, il faut cependant distinguer entre les deux qualités des racines que sont leurs natures (réelles, négatives) et leurs multiplicités. Dans le cas où les racines sont imaginaires ou réelles positives, la méthode s'applique et les oscillations, s'exprimant par une exponentielle réelle $e^{\rho t}$, ne sont pas bornées. En cas d'occurrence de racines multiples, la méthode ne permet plus l'obtention de *n* équations indépendantes et l'expression des solutions n'est plus valable. Lagrange considère cependant, reprenant un raisonnement sur les infiniment petits élaboré par d'Alembert pour le cas d'une corde chargée de deux masses, qu'en cas de racines multiples le « temps *t* sort du sinus » impliquant des solutions instables de la forme $t\sin(\rho t+\beta)$. Si ce type de solutions peut effectivement intervenir dans la résolution d'un système d'équations différentielles d'ordre *n* dont la matrice aurait des valeurs propres multiples, il n'intervient pas dans l'hypothèse symétrique ou plus généralement

De là on tire une méthode générale pour voir si l'état d'équilibre d'un système quelconque donné de corps est stable, c'est-à-dire si, les corps étant infiniment peu dérangés de cet état, ils y reviendront d'eux-mêmes, ou au moins tendront à y revenir. […]

1° Si toutes les racines de cette équation sont réelles négatives et inégales, l'état d'équilibre sera stable en général, quel que soit le dérangement initial du système.

2° Si ces racines sont toutes réelles, positives ou toutes imaginaires ou en partie positives, et en partie imaginaires, l'état d'équilibre n'aura aucune stabilité, et le système une fois dérangé de cet état ne pourra le reprendre ;

3° Enfin, si les racines sont en partie réelles négatives et inégales, et en partie réelles négatives et égales ou réelles et positives, ou imaginaires, l'état d'équilibre aura seulement une stabilité restreinte et conditionnelle. [Lagrange 1766, p. 532].

La « méthode générale » donnée par Lagrange pour caractériser la stabilité des systèmes mécaniques ne sera pas remise en cause avant la publication du mémoire de Weierstrass de 1858 intitulé « Ueber ein die homogenen functionen zweiten grades betreffendes Theorem, nebst anwendung desselben auf die theorie der kleinen schwingungen » et qui, comme nous l'avons vu, contient pour Kronecker un « germe » de la notion de diviseur élémentaire .

Après avoir indiqué et énoncé la forme des intégrales, Lagrange a conclu que, comme les oscillations $x_1$, $\frac{dx_1}{dt}$ restent toujours petites si elles le sont à l'origine, l'équation ne peut pas avoir de racines égales car les intégrales pourraient devenir arbitrairement grandes avec le temps. La même affirmation se trouve répétée chez Laplace lorsqu'il traite dans la *Mécanique céleste* des variations séculaires des planètes. Beaucoup d'autres auteurs, comme, par exemple, Poisson, mentionnent cette même conclusion. Mais cette conclusion n'est pas fondée […] on peut énoncer le même résultat, que les racines de l'équation $f(s)=0$ soient ou non toutes distinctes ; l'homogénéité de cette conclusion n'a pu être découverte dans la passé car on a toujours envisagé ce cas [des racines multiples] par des approches particulières. [Weierstrass 1858, p. 244, traduction F.B.]([41]).

Comme nous l'avons déjà mentionné, cette conclusion homogène est fondée sur l'énoncé d'un théorème portant sur les couples $(\Phi,\Psi)$ de fonctions homogènes du second degré et s'inscrivant dans la théorie des formes quadratiques ([42]), dans le cas où la forme $\Phi$ est définie, positive et où les coefficients de $\Phi$ et $\Psi$ sont réels, le théorème de Weierstrass énonce que les racines $s_i$ de l'équation $f(s)=0$ sont toujours réelles indépendamment de leur diagonalisable. L'absence de notion d'indépendance linéaire ne permet pas d'imaginer deux solutions s'exprimant à partir d'une même fonction $e^{\rho t}$ et tout de même « indépendantes ». Il y a ici deux familles indépendantes : celle des fonctions $(e^{\rho i t})$ dans l'espace des fonctions dérivables et celle des vecteurs de Àⁿ auquel est isomorphe l'espace des solutions de l'équation. C'est l'indépendance de cette dernière famille qui est nécessaire et qui peut donc être réalisée pour des racines multiples.

[41] Nachdem Lagrange die Form der Integral angegeben und gezeigt hat, wie die willkürlichen Constanten derselben durch die Anfangswerthe von $x_1$, $\frac{dx_1}{dt}$, u.s.w. stets unendlich klein bleiben, wenn sie es ursprünglich sind, auch die an, dass die genannte Gleichung keine gleiche Wurzeln haben dürfe, weil sonst in den Integralen Glieder vorkommen würden, die mit der Zeit beliebig gross werden könnten. Dieselbe Behauptung findet sich bei Laplace wiederholt, da wo er in der Mécanique céleste die Säcular-Störungen der Planeten behandelt, und ebenso, so viel mir bekannt ist, bei allen übrigen diesen Gegenstand behandelnden Autoren, wenn sie überhaupt den Fall der gleichen Wurzeln erwähnen, was z.B. bei Poisson nicht geschieht. Aber sie ist nicht begründet. [...], wenn nur die Function $\Psi$ stets negativ bleibt, und ihre Determinante nicht Null ist, was stattfinden kann, ohne dass die Wurzeln der Gleichung $f(s) = 0$ alle von einander verschieden sind ; wie man denn auch wirklich besondere Fälle der obige, Gleichungen, bei denen die Bedingung nicht erfüllt ist, mehrfach behandelt und doch keine Glieder von der angegebene Beschaffenheit gefunden hat.

[42] La forme $\Phi$ étant supposée de déterminant non nul. La relation entre les travaux de Weierstrass et l'arithmétique des formes quadratiques sera développée dans la suite de cet article. Voir [Hawkins 1977, p. 177] à propos de l'influence possible, sur la démonstration par Weierstrass de la stabilité, de [Dirichlet 1846] publié dans la troisième édition de la *Mécanique Analytique* [Lagrange 1853]. Dirichlet y propose une démonstration basée sur la définition de la continuité de Cauchy de la propriété selon laquelle un état d'équilibre d'un système mécanique conservatif est stable si la fonction potentielle atteint un maximum strict.

multiplicité. La stabilité des systèmes mécaniques est donc assurée par leur nature quadratique. Mais comme le montre l'appel d'Yvon-Villarceau qui, comme nous l'avons vu, attirera en 1870 l'attention des géomètres de l'Académie sur une « lacune » dans la relation établie par Lagrange entre stabilité mécanique et multiplicité des racines caractéristiques, le problème des petites oscillations est considéré comme résolu dans la période qui sépare Lagrange des théorèmes de Weierstrass (1858 et 1868) et de Jordan (1871 et 1872). Comment un problème considéré comme résolu est-il à l'origine de la discussion mathématique que forme notre corpus ?

Le problème des petites oscillations est d'abord apparu dans l'œuvre de Lagrange en 1766 pour l'étude d'un fil chargé d'un « nombre quelconque de masses ». La version de la *Mécanique analytique* dépasse cependant le cas des cordes vibrantes. Entre 1766 et 1788, l'efficacité de son procédé d'intégration des systèmes différentiels à coefficients constants conduit Lagrange à donner une même mathématisation au problème des petites oscillations d'un fil et au problème des petites oscillations des planètes sur leurs orbites :

> Si les Planètes étaient simplement attirées par le Soleil, et n'agissaient point les unes sur les autres, elles décriraient autour de cet astre, des ellipses variables suivant les lois de Kepler, comme Newton l'a démontré le premier, et une foule d'Auteurs après lui. Mais les observations ont prouvé que le mouvement elliptique des Planètes est sujet à de petites oscillations, et le calcul a démontré que leur attraction mutuelle peut en être la cause. Ces variations sont de deux espèces : les unes périodiques et qui ne dépendent que de la configuration des Planètes entre elles ; celles-ci sont les plus sensibles, et le calcul en a déjà été donné par différents Auteurs ; les autres séculaires et qui paraissent aller toujours en augmentant, ce sont les plus difficiles à déterminer tant par les observations que par la Théorie. Les premières ne dérangent point l'orbite primitive de la Planète ; ce ne sont, pour ainsi dire, que des écarts passagers qu'elle fait dans sa course régulière, et il suffit d'appliquer ces variations au lieu de la Planète calculé par les Tables ordinaires du mouvement elliptique. Il n'en est pas de même des variations séculaires. Ces dernières altèrent les éléments mêmes de l'orbite, c'est-à-dire la position et la dimension de l'ellipse décrite par la planète ; et quoique leur effet soit insensible dans un court espace de temps, il peut néanmoins devenir à la longue très considérable. [Lagrange 1781, p. 125].

Malgré les approximations nécessaires à la représentation des « équations séculaires des mouvements des nœuds et des inclinaisons des planètes » par « autant d'équations linéaires du premier ordre qu'il y a d'inconnues », la mathématisation élaborée par Lagrange est adoptée par Pierre-Simon Laplace dès 1775 en raison de l'efficacité de la « méthode fort ingénieuse » ramenant la résolution de problèmes mécaniques à la considération d'une équation algébrique et la question de la stabilité des inégalités séculaires à la nature des racines de cette équation. Dans deux publications successives de 1787 et 1789, Laplace entreprend de démontrer mathématiquement, en « dehors de toute hypothèse » numérique sur la masse des planètes, la stabilité du système solaire c'est-à-dire l'occurrence de racines réelles, négatives et inégales pour l'équation caractéristique de la méthode élaborée par Lagrange ([43]). La démonstration de 1789 repose sur un « artifice particulier » s'appuyant sur la présence de « rapports remarquables » [Laplace 1789, p. 297], $(1,2) = (2,1)$ et $[1,2] = [2,1]$ dans le système différentiel :

$$0 = (1)\,\frac{d^2\xi}{dt^2} + (1,2)\,\frac{d^2\psi}{dt^2} + (1,3)\,\frac{d^2\varphi}{dt^2} + \&c + [1]\,\xi + [1,2]\,\psi + [1,3]\varphi + \&c.$$

$$0 = (2)\,\frac{d^2\psi}{dt^2} + (1,2)\,\frac{d^2\xi}{dt^2} + (2,3)\,\frac{d^2\varphi}{dt^2} + \&c + [2]\psi + [1,2]\xi + [2,3]\varphi + \&c$$

[43] La preuve de Laplace, basée sur la conservation de l'intégrale première donnant l'énergie mécanique d'un système $L\text{-}V=C$, présuppose la forme erronée donnée par Lagrange aux solutions en cas de racines multiples (voir la note 37), elle est donc circulaire pour ce cas de figure. Pour une formulation de la démonstration de Laplace dans le cadre des mathématiques contemporaines, voir [Hawkins 1975, p. 15].

$$0 = (3)\ \frac{d^2\varphi}{dt^2} + (I,3)\ \frac{d^2\xi}{dt^2} + (2,3)\ \frac{d^2\psi}{dt^2} + \&c + [3]\varphi + [I,3]\xi + [2,3]\psi + \&c$$

Comme le manifestent les titres des travaux d'auteurs du corpus comme le mémoire de Cauchy de 1829 intitulé « Sur l'équation à l'aide de laquelle on détermine les inégalités séculaires des planètes », celui de Sylvester de 1852 « Sur une propriété nouvelle de l'équation qui sert à déterminer les inégalités séculaires des planètes » ou encore la formulation abrégée d'Hermite en 1857 « Mémoire sur l'équation à l'aide de laquelle, etc. », les travaux de mécanique de Lagrange et Laplace confèrent une identité spécifique à une équation algébrique qui se caractérise par la nature de ses racines et les « rapports remarquables » du système dont elle est issue. La spécificité de cette équation donne au corpus son principal élément de cohérence interne. Cauchy, par exemple, reconnaît une analogie formelle entre différents problèmes, comme les petites oscillations des systèmes mécaniques et la classification des coniques et quadriques relevant d'une recherche d'extremum d'une fonction homogène du second degré ([44]), en raison des caractéristiques « dignes de remarque » des équations que leurs résolutions mettent en œuvre [Cauchy 1829, p. 173]. Au cours de la période qui sépare les travaux de Lagrange de 1766 et les organisations théoriques élaborées par Jordan et Kronecker en 1874, l'« équation à l'aide de laquelle on détermine les inégalités séculaires des planètes » revêt une identité qui ne se limite pas à un cadre théorique mais l'identifie à un corpus d'auteurs que l'on sollicite lorsque cette équation se manifeste et qui s'étoffe en conséquence.

### 4. L'identité algébrique d'une pratique portée par la discussion.

Dans un article proposant une étude détaillée des principaux textes de notre corpus, nous avons montré que la spécificité de l'équation des petites oscillations est indissociable de celle d'une pratique consistant à exprimer les solutions d'un système linéaire par des factorisations polynomiales de l'équation caractéristique de ce système, c'est-à-dire par les expressions polynomiales (*) déjà mentionnées qui mettent en jeu les quotients $\frac{\frac{P_{1i}}{S}(x)}{x - s_j}$ ainsi que les quotients successifs $\frac{\Delta_{i+1}}{\Delta_i}$ des mineurs principaux de *S*. [Brechenmacher 200?]. Elaborée à l'origine par Lagrange comme spécifique aux problèmes de petites oscillations, cette pratique originale s'affranchit de la méthode des coefficients indéterminées et, en un jeu sur les primes et les indices des coefficients des systèmes linéaires est à l'origine d'une caractéristique des systèmes issus de la mécanique, les « rapports remarquables » de la disposition en miroirs de systèmes que nous désignerions aujourd'hui comme symétriques. Le corpus de la discussion sur les racines de l'équation des petites oscillations manifeste un héritage de la pratique élaborée par Lagrange sur une période antérieure à l'élaboration des théories algébriques qui donneront à cette pratique l'identité d'une méthode de transformation d'un système linéaire par la décomposition de la forme polynomiale de l'équation caractéristique associée ([45]). Caractérisée par une pratique

[44] Si $q$ est une forme quadratique et $\varphi$ sa forme polaire alors $dq_x(y)=2\varphi(x,y)$, on peut déterminer un vecteur propre $e$ de $q$ par une recherche d'extremum sur la sphère unité.

[45] Comme nous l'avons vu dans le paragraphe II. 1. et la note n°23, on interprèterait aujourd'hui la pratique de Lagrange comme une méthode revenant à donner une expression polynomiale générale des vecteurs propres d'une matrice comme quotients de mineurs extraits du déterminant caractéristique $|A-\lambda I|$ et d'une factorisation de l'équation caractéristique par un terme linéaire : $|A-\lambda I|/(\lambda-\lambda_i)$ (où $\lambda_i$ est une racine caractéristique de $A$). Une telle expression polynomiale s'identifie à un facteur près aux colonnes non nulles de la matrice des cofacteurs de $|A-\lambda I|$. Une telle formulation introduit cependant dans le discours historique des aspects propres à la théorie des matrices, comme les termes « transformation » et « symétrique », qui sont absents du texte de Lagrange et masquent le caractère spécifique d'une pratique originale. Lorsque Lagrange ramène l'intégration d'un système de n équations à

indissociable d'une équation, l'identité de notre corpus présente donc une nature algébrique que nous allons préciser.

A l'exception de Jordan et Darboux, aucun des auteurs du corpus ne se réfère à une théorie algébrique. L'identité algébrique de la discussion ne peut donc se caractériser, avant 1874, comme une identité théorique, elle se place davantage du côté des pratiques que des méthodes. A la permanence d'une pratique polynomiale de résolution des systèmes linéaires il faut en effet opposer la diversité des méthodes élaborées dans les cadres théoriques distincts au travers desquels se déploie la discussion. La méthode élaborée par Lagrange pour ramener l'intégration d'un système à celle de *n* équations indépendantes est indissociable d'une représentation mécanique envisageant les petites oscillations d'un système comme une combinaison d'oscillations indépendantes, cette méthode s'avère très différente de celle de Cauchy, sous tendue par des représentations géométriques de changements d'axes orthogonaux. Dans le cadre du corpus, une même pratique algébrique est employée aux seins de méthodes et cadres théoriques différents dans lesquels elle adopte des représentations diverses. Le point de vue selon lequel les procédés algébriques ne portent pas de représentations propres mais prennent leurs « significations » relativement aux théories dans lesquels ils sont employés est, comme nous l'avons vu, celui que Kronecker oppose en 1874 à l'organisation algébrique que veut donner Jordan à la théorie des formes bilinéaires. Chez Kronecker comme chez la plupart des auteurs du corpus, le travail algébrique n'est pas une fin en soi mais sert des objectifs issus de contextes divers et de disciplines comme la mécanique, la géométrie ou l'arithmétique. Si, parallèlement à cette variabilité des contextes, une identité algébrique se manifeste dans la discussion sur l'équation des petites oscillations, celle-ci ne tient pas à la perception d'une théorie algébrique sous jacente mais à une perspective historique. C'est en effet par des références à des travaux plus anciens que les auteurs du corpus comme Laplace, Cauchy ou Weierstrass donnent à leurs travaux une identité dépassant les contextes dans lesquels ceux-ci s'inscrivent. Cauchy par exemple, dépasse en 1829 le contexte géométrique de la recherche des axes principaux des coniques et quadriques par sa référence à Lagrange qui identifie la pratique algébrique dont il s'inspire. Avant que le corpus ne soit identifié à l'histoire d'une théorie (comme Kronecker qui se réfère à la théorie des faisceaux de formes quadratiques), c'est une identité historique qui permet aux auteurs de considérer leurs travaux comme relevant d'une discussion de nature algébrique.

De l'origine de la discussion chez Lagrange à ses deux fins chez Weierstrass et Jordan, la revendication de généralité est la principale motivation qui amène des auteurs à insérer leurs propres travaux dans une perspective historique plus large. Dès 1766, Lagrange élabore la pratique à l'origine de la discussion dans une ambition de généralisation du traitement par Jean le Rond d'Alembert des petites oscillations d'un fil chargé de deux ou trois masses. La capacité de sa méthode à résoudre un problème mécanique pour un « système quelconque » de corps lui confère une portée générale qui justifie la place donnée au problème comme premier exemple d'application de la *Mécanique analytique* de 1788. C'est également la généralité du problème qui suscite une discussion qualitative sur la nature des racines d'une équation de degré *n* ne pouvant être résolue de manière effective. La perception de ce qui s'avère « général » dans cette discussion évolue avec le temps. Pour

---

celle de n équations indépendantes, le changement de forme du système n'est pas envisagé comme le résultat une transformation mais s'avère indissociable d'une représentation mécanique implicite selon laquelle les petites oscillations d'un système de *n* corps se comportent comme si elles étaient composées de *n* mouvements simples. Ce que nous désignerions aujourd'hui comme le caractère symétrique des systèmes linéaires mécaniques se manifeste comme une conséquence de la pratique algébrique originale élaborée par Lagrange et consistant, en un jeu sur les primes et les indices des coefficients des systèmes linéaires que l'on interpréterait aujourd'hui dans le cadre d'une orthogonalité duale, à réduire l'intégration de ces derniers à la décomposition algébrique d'une équation. Pour une description détaillée de la pratique de Lagrange et de son héritage chez Cauchy, Hermite et Weierstrass voir [Brechenmacher 200?].

Lagrange, la nécessité d'obtenir des racines caractéristiques distinctes ne réduit pas la généralité polynomiale car ce cas de figure est considéré comme mécaniquement incorrect puisque impliquant des oscillations croissantes à l'infini. Les présupposés de stabilité mécanique impliquent des propriétés des racines de l'équation associée- les oscillations sont petites donc les racines sont réelles, négatives et inégales -, les enjeux de cette implication changent cependant de nature lorsque Lagrange aborde, en 1774, le problème des variations séculaires des planètes. Assurer la stabilité du système du monde suscite la recherche par Laplace d'une démonstration générale de la réalité des racines. L'intervention de Cauchy dans la discussion est d'abord motivée par la généralisation à *n* variables, que permet le problème des petites oscillations, d'une méthode élaborée pour deux ou trois variables dans un cadre géométrique. Les « formules algébriques » généralisant les expressions de changements d'axes orthonormés ne « subsistent » cependant pas à la « condition » de multiplicité des racines caractéristiques conduisant à la forme $\frac{0}{0}$. Pour Cauchy, qui critique les méthodes qui tendent « à faire attribuer aux formules algébriques une étendue indéfinie, tandis que […] la plupart de ces formules subsistent uniquement sous certaines conditions », le problème posé par l'occurrence de racines multiples manifeste un cas singulier limitant l'« étendue » d'une formule et nécessitant un traitement spécifique initialement mené par introduction d'infiniment petits. C'est par opposition à de tels cas singuliers qui « encombrent » les formules de l'algèbre que Cauchy introduit dès 1826 le calcul des résidus à l'aide duquel il donnera en 1839 une résolution homogène au problème des petites oscillations. Le renouvellement de la signification de la généralité qui accompagne cette exigence d'homogénéité va porter la discussion des petites oscillations de Cauchy à Weierstrass, en passant par Jacobi, Borchardt, Hermite, Johann Dirichlet et Sylvester ([46]). La question se focalise alors sur le problème de la multiplicité des racines, d'abord abordé par des méthodes différentielles par Sturm, Cauchy et Sylvester. Ce dernier, après une première série de travaux consacrée à la nature des racines algébriques (1840-1842) articule la méthode différentielle à la théorie du déterminant et ses représentations géométriques en termes de « rectangle symétrique » par rapport à sa « diagonale ». Ses travaux qui, avec ceux de Cayley, sont à l'origine de la théorie des invariants et voient l'introduction des termes « matrices » et « mineurs » (1850-1852) ([47]), sont investis dans un cadre arithmétique en relation avec les préoccupations d'Hermite sur les formes quadratiques et la décomposition en quatre carrés (1853-1856).

La question de la généralité est indissociable de l'identité algébrique du corpus mais tout en donnant, comme nous l'avons vu, une conclusion algébrique, générale et homogène à la discussion, Weierstrass s'appuie en 1858 sur des procédés de « transformation » de « formes » linéaires qui manifestent l'orientation arithmétique qu'ont pris les travaux sur l'équation des petites oscillations dans les années 1850. Dans un cadre arithmétique, le problème est présenté sous une forme semblable à la loi d'inertie, « réduire un polynôme homogène du second degré à une somme de carrés », et abordé par des méthodes propres à la théorie des formes quadratiques. On recherche les invariants permettant de caractériser les classes d'équivalences des formes. Comme le formule clairement Darboux qui, en 1874, généralise des procédés élaborés par Hermite pour donner une nouvelle démonstration du théorème de Weierstrass de 1858, l'invariant du couple de formes $(A(x_1,\ldots,x_n), x_1^2+\ldots+x_n^2)$ qu'est le déterminant $S(\lambda)=|A(x_1, x_2,\ldots,x_n)-\lambda(x_1^2+x_2^2+\ldots+x_n^2)|$ ne permet pas de caractériser les différentes classes d'équivalences et il faut aller au-delà du déterminant en considérant

[46] Le calcul des résidus de Cauchy permet notamment d'explorer les pôles de la fonction rationnelle $P_{ij}(x)/S(x)$ (où $S(x) = |A+xB|$ et $P_{ij}$ son mineur obtenu en développant la $i^e$ ligne et $j^e$ colonne), démarche à la base du théorème des diviseurs élémentaires. Voir à ce sujet [Hawkins, 1977, p.130].

[47] Au sujet des origines de la théorie des invariants, consulter [Parshall 1989 et 2006].

ses « mineurs » ([48]) [Darboux 1874, p. 367]. Contrairement à la loi d'inertie de la théorie arithmétique des formes quadratiques, le problème porte ici sur la transformation simultanée d'un couple de formes envisagé comme un polynôme de formes quadratiques $A-\lambda I$ et, en raison de ce caractère polynomial, s'insère pour Darboux dans « la théorie algébrique des formes quadratiques ». Si, dans un cadre arithmétique, le terme « forme » bénéficie d'une définition mathématique explicite en relation avec la notion d'équivalence, les significations des termes « formes » et « transformations » restent le plus souvent implicites et évoluent considérablement dans le cadre algébrique de la discussion des petites oscillations. Dans la *Mécanique analytique,* l'existence d'une forme intégrable du système des petites oscillations est garantie par la représentation mécanique du problème selon laquelle le mouvement d'une corde chargée de *n* masses se décompose mécaniquement en oscillations de *n* cordes chargées d'une seule masse. Si la méthode de Lagrange consiste bien à ramener l'intégration d'un système d'équations différentielles à celle de *n* équations indépendantes, le changement de la forme du système ne s'appuie pas sur un procédé de transformation mais sur la résolution d'une équation algébrique donnant les paramètres mécaniques du système. Chez Cauchy la notion de « transformation d'une fonction homogène » est indissociable d'une représentation géométrique de transformation des axes de coordonnées d'une conique. Si la question qui semble aujourd'hui essentielle, et que l'on formulerait comme la caractérisation de la forme d'une matrice non nécessairement diagonalisable car présentant des valeurs propres multiples, n'est pas abordée par Lagrange, Laplace ou Cauchy il ne s'agit pas là d'un défaut de généralité comme le sanctionne Kronecker en 1874 mais d'une question étrangère à la manière dont sont envisagés les termes « formes » et « transformations » d'un système. Cette question est au contraire naturelle dans le cadre arithmétique des travaux d'Hermite et Weierstrass des années 1850 qui ouvrent la voie aux premiers travaux de Christoffel sur les formes bilinéaires comme dans le cadre algébrique des travaux de Jordan ([49]). Dans la résolution générale que donne ce dernier au problème de l'intégration des systèmes d'équations linéaires à coefficients constants, la décomposition en facteurs irréductibles du polynôme caractéristique permet de regrouper les variables en « un certain nombre de séries » et de « ramener » le système d'équations à une « suite » de « formes simples » dont l'intégration est connue [Jordan 1871, 787]. Aux diverses significations associées au terme « formes » chez Lagrange, Laplace ou Cauchy, répond chez Hermite, Weierstrass, Kronecker, Jordan ou Darboux une mathématique dont la « forme », la « transformation » est l'objet. Cette mathématique qui permet de poser, en toute généralité, la question de la caractérisation de la forme d'un système linéaire est-elle une algèbre ou une arithmétique ? Comme nous allons le voir, la rencontre des héritages portés par la théorie arithmétique des formes quadratiques et la pratique algébrique spécifique à la discussion sur l'équation à l'aide de laquelle on détermine les inégalités séculaires des planètes suscite une tension entre arithmétique et algèbre. Si Jordan et Kronecker se réfèrent, lors de la querelle de 1874, à une histoire commune dans laquelle ils identifient une pratique spécifique, consistant à aborder les transformations des couples de formes bilinéaires *(A,B)* par la décomposition de la forme polynomiale de l'équation $|A+sB|=0$, les deux géomètres opposent la généralité de deux méthodes (calculs d'invariants, réductions canoniques), deux cadres théoriques (formes bilinéaires, substitutions), deux disciplines (arithmétique, algèbre) dans lesquels s'insèrent cette pratique partagée par les deux savants.

[48] Darboux voit notamment la notion de diviseur élémentaire de Weierstrass en germe dans les invariants introduits par Sylvester en 1851 pour caractériser les intersections de coniques et quadriques au sein de travaux qui voient la définition des « mineurs » d'une « matrice ». Voir à ce sujet [Brechenmacher 2006d].

[49] Dans deux mémoires successifs, consacrés à ce que nous désignerions comme le cas hermitien dans lequel les racines restent réelles, et publiés en 1864, Christoffel s'appuie sur les résultats de Weierstrass de 1858 pour généraliser des travaux mécaniques de Clebsch [1860] et des résultats d'[Hermite, 1854, 1855, 1856] sur les résidus biquadratiques, les propriétés arithmétiques des nombres $a+bi$ ($a$ et $b$ entiers) et la décomposition d'un entier en somme de quatre carrés. Le premier mémoire se présente comme un développement mathématique du second, consacré à des problèmes de petites oscillations dans le cadre de la théorie de Cauchy selon laquelle la lumière correspond aux petites vibrations de molécules ponctuelles d'éther soumises à des forces attractives et répulsives. Le procédé nécessite l'intégration d'un système d'équations différentielles linéaires dont les coefficients sont des constantes complexes [Mawhin, 1981] présentant des « propriétés remarquables car ils apparaissent comme des extensions de ces équations qui interviennent dans la théorie des perturbations séculaires des planètes et dans tant d'autres recherches » [Clebsch 1860, p.326] et généralisant le cas quadratique à ce que Christoffel décrit comme des « fonctions bilinéaires » en référence à un mémoire de 1857 dans lequel Jacobi généralisait la loi d'inertie des formes quadratiques aux « fonctions bilinéaires » $\Sigma a_{ij}x_iy_j$.

## III. Les idéaux disciplinaires opposés par la querelle.

### 1. L'arithmétique des formes de Kronecker : homogénéité et effectivité.

L'origine des travaux sur les formes bilinéaires à Berlin dans les années 1860-1870 se caractérisait comme une entreprise théorique guidée par un idéal de généralité [Brechenmacher 2006a, p. 205], la résolution complète de la question de la transformation des couples de formes bilinéaires et quadratiques par les publications conjointes de Weierstrass et Kronecker de 1868 était un des fondements de la théorie. Pourquoi, dans ce cas, Kronecker publie-t-il en 1874 de nouveaux mémoires, objets de communications mensuelles à l'Académie durant l'hiver, sur une question qu'il considère lui-même complètement résolue depuis 1868 ? Kronecker présente sa nouvelle approche comme provenant d'échanges avec Ernst Kummer sur la nécessité de rassembler les différents résultats obtenus dans les années 1860 en une théorie arithmétique des « faisceaux de formes quadratiques » ([50]). Le qualificatif « arithmétique » s'associe à la revendication d'un héritage de l'arithmétique des formes quadratiques de Gauss et des travaux conjoints d'Hermite et de Kronecker au début des années soixante sur la résolution des équations algébriques par les équations modulaires des fonctions elliptiques. Cet héritage se manifeste notamment dans l'emploi du terme « forme » pour ce que d'autres désignent comme des « fonctions » ([Weierstrass 1858], [Christoffel 1866]) ou des « polynômes » ([Jordan 1873]). Le terme « forme », tel qu'introduit par Gauss dans les *Disquitiones arithmeticae* de 1801 ([51]), est associé à la notion arithmétique de classe d'équivalence dont Kronecker réalise la généralisation aux faisceaux de formes comme une « application des notions de l'arithmétique à l'algèbre » :

> En appliquant les notions de l'Arithmétique à l'Algèbre, on peut appeler *équivalentes* deux formes bilinéaires, dont l'une peut être transformée en l'autre par une même substitution, opérée sur les deux systèmes de variables, et ensuite on peut réunir en une même *classe* toutes les formes équivalentes.
>
> Pour que deux formes bilinéaires $\varphi(x,y)$ et $\psi(x,y)$ appartiennent à une même classe, il faut et il suffit que les deux faisceaux formés des deux paires de fonctions conjuguées $u\varphi(x,y)+v\varphi(y,x)$, $u\psi(x,y)+v\psi(y,x)$ soient équivalents. [Kronecker 1874b, p. 415]

Les résultats sont exposés pour les formes quadratiques mais les méthodes peuvent s'appliquer au « cas particulier » des formes bilinéaires par la distinction de deux types de relations selon que l'on fasse opérer des substitutions linéaires ou orthogonales, ce sont les relations d' « équivalence » et de « congruence » des formes ([52]). Le problème essentiel de la théorie se formule désormais comme une question arithmétique de caractérisation de classes d'équivalences, cette nouvelle formulation s'accompagne d'idéaux disciplinaires qui se manifestent dans la critique de la forme canonique de Jordan :

> Depuis ma présentation « sur les formes quadratiques et bilinéaires » du 16 février, deux publication de M. C. Jordan sont parues sur le même objet […] Sa proposition selon laquelle « la condition suffisante et nécessaire pour l'équivalence du système de deux formes est l'identité des réduites », bien que parfaitement exacte, n'est pas d'un contenu satisfaisant. Elle n'énonce en effet aucun procédé pratique pour caractériser l'équivalence d'un système de formes et doit être distinguée de la possibilité immédiate qu'offre le critère théorique d'équivalence, de former un système complet d'invariants –au sens propre du mot-à partir des coefficients des formes données. [Kronecker 1874c, p. 382, traduction F.B.] ([53])

[50] Kronecker démontre ainsi que le problème I de Jordan est un cas particulier de l'étude d'un faisceau $u\varphi+v\psi$ où $\varphi$ est bilinéaire et $\psi$ est quadratique (la forme $\psi$ est par exemple, pour le cas de la classification des coniques, la forme identité qui correspond à la condition d'orthonormalité des substitutions). C'est essentiellement pour cette raison que Kronecker rejette la classification de Jordan en trois problèmes, qui s'appuie sur l'action des groupes classiques.

[51] Consulter à ce sujet la récente thèse de doctorat de M. Bullynck, [Bullynck 2006].

Dans cet extrait, Kronecker formule tout à la fois sa critique principale de Jordan et l'ambition de ses travaux de 1874. Son approche arithmétique s'accompagne d'un idéal d'effectivité condamnant la réduction canonique qui nécessite l'extraction des racines d'une équation algébrique générale, l'équation caractéristique, et ne permet donc de procédé « pratique » que pour les équations de degré inférieur à 5. Cette critique atteint aussi bien la forme de Jordan que les diviseurs élémentaires tels qu'introduits par Weierstrass en 1868 par des factorisations du polynôme caractéristique en expressions linéaires et dont Kronecker donne une nouvelle définition effective ([54]):

> Dans la théorie arithmétique des formes, il faut à vrai dire se contenter de l'indication d'une méthode pour décider de la question de l'équivalence […] (cf. Gauss : Disquitiones arithmeticae, Sectio V [...]). La méthode elle-même peut nécessiter la transition à des formes réduites : cependant il ne faut pas omettre de voir que celles-ci ont une toute autre importance dans la théorie arithmétique des formes qu'en algèbre. En effet, comme les invariants des formes équivalentes ne sont, de part leur natures que des fonctions arithmétiques des coefficients, il ne faut pas s'étonner que ceux-ci sont certes définis directement sans l'être explicitement, mais ne pourront être représentés que comme des résultats d'opérations arithmétiques; la plupart des notions arithmétiques se comporte de manière tout à fait semblable, comme la notion simple de plus grand diviseur commun. [*Ibidem*] ([55]).

---

[52] Cette dénomination ne revêt pas un sens identique à celui des termes « équivalence » et « congruence » dans les mathématiques d'aujourd'hui. De telles définitions seront données quelques années plus tard par [Frobenius 1878]. Kronecker appelle ainsi « équivalents » aussi bien des formes congruentes *P* et *P'* que des polynômes de formes équivalents *pP+qQ* et *pP'+qQ'*. Dans la même année, il introduira le terme « formes correspondantes » pour le premier cas [Kronecker 1874d, p. 423]. La dénomination de formes « congruentes » sera réservée à des formes « correspondantes » par l'action du groupe orthogonal.

[53] Seit meinem am 16. Februar gehaltenen Vortrage « über quadratische und bilineare Formen » sind zwei Publikationen des Hrn. C. Jordan über den selben Gegenstand erschienen [....].Seine Proposition « dass für die Aequivalenz der Systeme zweier Formen die Uebereinstimmung der Reductirten nothwendig und hinreichend sei », ist zwar vollkommen richtig, aber zu dürftigen Inhalts, denn es handelt sich nicht um die Angabe eines praktischen Verfahrens zur Entscheidung der Frage der Aequivalenz gegebener Formensysteme, sondern um eine möglichst unmittelbare Anknüpfung der theoretischen Kriterien der Aequivalenz an die Coëfficienten der gegebenen Formen, d.h. die Aufstellung eines vollständigen Systems von « Invarianten », im höheren Sinne des Wortes.

[54] Ce passage est extrait d'une note de bas de page appelée, dans le texte original, par un astérisque *) placé à la fin de la citation précédente.

[55] In der arithmetischen Theorie der Formen muss man sich freilich mit der Angabe eines Verfahrens zur Entscheidung der Frage der Aequivalenz begnügen und das betreffende Problem wird deshalb auch ausdrücklich in dieser Weise formuliert (cf. Gauss : Disquitiones arithmeticae, Sectio V [...]) Das Verfahren selbst beruht auch dort auf dem Uebergange zu reductiren Formen : doch ist dabei nicht zu übersehen, dass denselben in den arithmetischen Theorien eine ganz andere Bedeutung zukommt als in der Algebra. Da nämlich die Invarianten äquivalenter Formen dort ihrer Natur nach nur zahlentheoretische Functionen der Coëfficienten sind, so kann es nicht befremden, wenn dieselben zwar direct definiert aber nicht explicite sondern nur als Endresultate

Les invariants qui caractérisent les classes d'équivalences doivent s'exprimer de manière effective par des opérations arithmétiques sur les coefficients, ils ne peuvent pas être « représentés explicitement » mais sont définis par un procédé de calcul : ce sont les plus grands diviseurs communs des mineurs respectifs du déterminant caractéristique ([56]). L'idéal d'effectivité de Kronecker s'oppose à l'énoncé d'une « écriture littérale » comme la forme canonique de Jordan car la décomposition algébrique des polynômes dépend du « domaine de rationalité » sur lequel on travaille, au contraire de l'algorithme arithmétique du p.g.c.d.:

> Je remarque à cette occasion, que les invariants algébriques sont déduits de manière effective et dans leur pleine généralité comme plus grand communs diviseurs de fonctions entières et nullement, comme on l'a accepté jusqu'ici, par des écritures littérales. Je suis arrivé à ce résultat depuis quelques années par mon travail sur les discriminants des équations algébriques et plus tard, par mon travail sur les transformations linéaires soumis à l'Académie en octobre 1868. [Kronecker 1874a, p. 353, traduction F.B.] ([57]).

La querelle de 1874 voit Kronecker expliciter pour la première fois l'idéal d'effectivité souvent associé par l'historiographie à la théorie arithmétique des grandeurs algébriques des années 1880-1890. En 1874, la pratique des calculs de Kronecker ne reflète cependant pas la radicalité des discours. Si, comme nous l'avons vu dans la deuxième partie, un procédé « formel » ne doit en aucun cas participer de l'organisation théorique, dans la pratique de ses calculs et démonstrations, Kronecker recourt à des « formes normales » similaires aux formes canoniques de Jordan et nécessitant des extractions de racines d'équations algébriques. Le théorème principal de Kronecker caractérise les couples de formes bilinéaires ou quadratiques, singuliers ou non, par des « formes élémentaires » associées à des invariants. Chaque classe d'équivalence des faisceaux de formes bilinéaires est associée à une suite de diviseurs élémentaires, polynômes homogènes en deux variables *f(u,v)*. Le premier terme de la suite est le déterminant caractéristique du faisceau, le membre suivant est le « plus grand commun diviseur des premiers sous déterminants » et ainsi de suite en considérant les p.g.c.d. des sous déterminants successifs. Ces suites de polynômes se décomposent en des suites « élémentaires », correspondantes aux diviseurs élémentaires de Weierstrass, chaque faisceau de formes s'écrit alors comme un « agrégat » de « faisceaux élémentaires » donnant la structure des classes d'équivalences arithmétiques en un énoncé qui rassemble les résultats obtenus par Weierstrass et Kronecker en 1868 pour les cas singuliers et non singuliers ([58]) :

I. $$(-1)^n \sum_h x_h y_{h+1} + \sum_h (-1)^h y_h x_{h+1} + x_n y_n \qquad (h=0,1,\dots,n-1)$$

II. $$(-1)^m \sum_h x_h y_{h+1} + \sum_h (-1)^h y_h x_{h+1} \qquad (h=0,1,\dots,2m-2)$$

---

arithmetischer Operationen dargestellt werden können ; denn ganz ähnlich verhält es sich mit den meisten arithmetischer Begriffen, z.B. schon mit jenem einfachsten Begriffe des grössten gemeinsamen Theilers.

[56] On reconnaitrait aujourd'hui les facteurs invariants d'une matrice *A* sur un anneau principal. Si $D(\lambda)=|A-\lambda I|$, dans la suite des p.g.c.d. des mineurs successifs $\square_r(\lambda)$, $\square_{r-1}(\lambda)$, ..., $\square_1(\lambda)$, chaque polynôme est divisible par le précédent et si les quotients correspondants sont désignés par $i_1(\lambda)$, $i_2(\lambda)$,...$i_r(\lambda)$ et appelés les polynômes invariants de la matrice $A(\lambda)$ la décomposition de ces polynômes en facteur irréductibles distincts sur *K* donne les diviseurs élémentaires de la matrice $A(\lambda)$ sur le corps *K* (voir annexe 3).

$$i_1(\lambda)=[\varphi_1(\lambda)]^{s_1}\dots[\varphi_s(\lambda)]^{c_s},\ i_2(\lambda)=[\varphi_1(\lambda)]^{d_1}\dots[\varphi_s(\lambda)]^{d_s},\ \dots,\ i_r(\lambda)=[\varphi_1(\lambda)]^{l_1}\dots[\varphi_s(\lambda)]^{l_s}\ (c_k>d_k>\dots>l_k)$$

[57] Ich bemerke bei dieser Gelegenheit, dass, wie hier, so die algebraischen Invarianten überhaupt in ihrer wahren Allgemeinheit nur aus grössten gemeinsamen Theilern von ganzen Functionen gegebener Elemente herzuleiten und keineswegs, wie bisher angenommen wurde, durch literale Bildungen zu erschöpfen sind. Ich bin hierauf schon vor einer langen Reihe von Jahren bei meinen Untersuchungen über die Discriminante von algebraischen Gleichungen geführt worden, sowie später bei meiner Arbeit über lineare Transformationen, welche ich im October 1868 der Akademie mitgetheit habe.

[58] Pour une formulation de ce théorème du point de vue de la théorie des matrices des années 1930, voir [Aitken et Turnbull 1933, p. 120].

III $$a\sum_h x_h y_{h+1} + b\sum_h y_h x_{h+1} \quad (h=0,1,\dots,n-1\ ;\ a^2 \gtrless b^2)$$

Ces trois fonctions ne sont plus décomposables d'une manière analogue et c'est pourquoi je les désigne comme *formes élémentaires*. […] si l'on désigne par $F_1$, $F_2$, $F_3$ respectivement les faisceaux qui appartiennent aux formes élémentaires I, II, III et par $D_1$, $D_2$, $D_3$ leurs *déterminants,* on a

$D_1 = [u+ (-1)^n v]^{n+1}$,

$D_2 = [u+ (-1)^m v]^{2m}$,

$D_3 = (au+bv)^m(av+bu)^m$, ($n+1$ étant égal à $2m$)

$D_3 = 0$ ($n+1$ étant un nombre impair). [Kronecker 1874b, p. 418]

## 2. La réduction algébrique de Jordan : généralité et simplicité.

A la généralité, et aussi l'uniformité, revendiquées par Kronecker, Jordan oppose une constante revendication de simplicité qui, comme le manifestent les termes soulignés de la citation suivante, est indissociable d'un procédé de réduction :

Nous pensons donc satisfaire les géomètres en exposant, pour la solution de ces questions [les trois questions distinguées par Jordan dans la citation donnée dans la première partie], une méthode nouvelle très **simple**, et ne comportant plus aucun cas d'exception […]. La méthode par laquelle nous traitons le premier problème est plus **simple** que celle dont s'est servie M. Kronecker ; mais elle repose sur les mêmes principes, et la solution des deux autres questions s'en déduit fort **aisément** ; aussi aurions nous hésité à publier ce travail ; mais nos scrupules ont été levés par la lecture d'un Mémoire récent de M. Kronecker (*Monatsbericht*, mars 1874). Cet habile analyste, après avoir contesté ce que nous avions dit […] sur la **simplicité** des principes dont cette question dépend, annonce en effet, qu'il s'est occupé de la transformation d'un système en lui-même, mais qu'il n'a pas réussi à obtenir une conclusion complètement satisfaisante, et que, dans l'étude de ce problème, il n'a pu tirer que peu de profit de son procédé de **réduction**. Pour qu'un géomètre aussi exercé ait pu méconnaître ainsi, du même coup, la **simplicité** et la portée de sa propre méthode, il faut évidemment que la question ne soit pas encore suffisamment élucidée. […] Nous ferons observer en outre qu'on gagne beaucoup en **simplicité** et en **élégance** en opérant symétriquement sur les formes à **réduire**. Tout le raisonnement tient en quelques lignes. On obtient directement les faisceaux **élémentaires** tandis que M. Kronecker a besoin dans certains cas d'une **réduction** ultérieure. Enfin, on reconnaît, chemin faisant, de la manière la plus **simple**, que la forme des **réduites** est complètement déterminée […]. [Jordan 1874b, p. 13].

Comme nous allons le voir dans ce paragraphe, les arguments de Jordan ne se réduisent pas au simplisme naïf caricaturé par Kronecker mais manifestent un idéal indissociable d'une pratique de « réduction » d'un problème général en une chaîne de sous problèmes dont les maillons doivent être les « plus simples » comme dans les deux résultats essentiels que sont le théorème de décomposition des groupes (Jordan-Hölder) et le théorème de réduction à une forme canonique. Cette pratique se manifeste notamment dans la critique que fait Jordan de la réduction obtenue par Kronecker en 1868 pour les couples singuliers de forme. Cette dernière ne satisferait pas l'exigence de simplicité des « vraies réduites » dont toute réduction ultérieure doit être impossible :

Le lecteur nous sera gré de reproduite ici le résultat publié à cette époque par M. Kronecker :

« *Si deux formes quadratiques P et Q, à n variables, satisfont à la condition (P,Q)=0, on pourra les réduire toutes les deux, par un changement de variables à la forme*

*(1)* $f_1x_{m+1}+f_2x_{m+2}+\dots+f_mx_{2m}+\mathcal{F}$

(2) *𝒥 étant une forme quadratique des n-2 m-1 dernière variables, et* $f_1,\ldots f_m$ *étant des fonctions linéaires quelconques de toutes les variables. »*

[…] [Les réduites de Kronecker] contiennent encore des coefficients indéterminées qu'une réduction ultérieure doit faire disparaître. Elles ne peuvent servir à constater l'équivalence de deux systèmes de deux formes *P* et *Q*, $P_1$ et $Q_1$ ; car on peut trouver pour *P* et *Q*, d'une infinité de manières, une infinité d'expressions différentes de l'espèce […]. Enfin les expressions (1) ne mettent pas en évidence le caractère fondamental des vraies réduites d'être décomposables en fonctions partielles ne contenant chacune qu'une portion des variables. [*Ibidem*].

Paradoxalement, alors que Jordan revendique constamment la simplicité de ses méthodes et malgré le succès rencontré par la publication du *Traité des substitutions*, félicité par les plus grands savants européens et auquel l'Académie attribue le prix Poncelet en 1873 ([59]), la scène parisienne résonne de « bruits défavorables » sur les travaux mathématiques de Jordan, critiqués comme « inintelligibles » et comme n'ayant « sans doute pas la portée qu'on leur attribue ». En témoignent les inquiétudes nourries par Jordan pour sa candidature de 1875 à la section de géométrie et qui font l'objet d'une lettre adressée à Hermite le 2 décembre 1874. La controverse qui s'achève alors semble jeter une ombre sur la carrière de Jordan, elle constitue certainement une des origines de ses relations tendues avec Hermite qui qualifie l'étude de ses travaux de « tellement difficile et tellement pénible », refuse de soutenir sa candidature à l'Académie et rejette sa demande de traduire les mémoires de Weierstrass de 1858 et 1868 ([60]). Des travaux incompris, des méthodes contestées, des démonstrations du *Traité* remises en cause (Jordan est tout prêt d'entamer une seconde querelle avec Eugen Netto, un élève de Kronecker, dont un mémoire remet en cause la démonstration de Jordan du théorème de décomposition des groupes [Netto 1874]) ([61]), la situation de Jordan est critique. Par contraste, au début des années 1880 Jordan sera devenu l'un des patrons des mathématiques parisiennes ; dans l'intervalle, et par opposition à l'abstraction que ses contemporains associent à ses recherches sur les substitutions, Jordan aura diversifié ses travaux à des champs d'applications, souvent par emploi de la forme canonique des substitutions.

La querelle de 1874 intervient à une charnière de la carrière de Jordan qui, dans les années 1870, se détache progressivement de sa profession d'ingénieur des mines pour venir occuper des positions institutionnelles clés des mathématiques parisiennes. Cette évolution de carrière est indissociable d'une évolution profonde des travaux du savant. Après une première période consacrée à la théorie des substitutions et à quelques autres domaines comme la géométrie algébrique ou la topologie, les années 1870-1880 voient les recherches de Jordan prendre une nouvelle envergure par une diversification qui se nourrit de la capacité du savant à appliquer les notions et méthodes développées pour la théorie des groupes à des domaines variés comme les systèmes d'équations différentielles linéaires (1871), la théorie des formes bilinéaires et quadratiques (1872-1875), l'intégration algébrique des équations différentielles (1875-1878) et la théorie des nombres (1878-1907). Le regard critique de Kronecker révèle les idéaux implicites qui accompagnent dans les applications les notions et méthodes de la théorie des groupes. La « simplicité » en est le terme clef comme nous allons le voir en portant un regard rétrospectif sur les recherches sur les substitutions des années 1860 qui permettront à Jordan d'être considéré par ses successeurs comme un « grand algébriste » [Picard 1922, p. VIII], un précurseur incompris en son temps et isolé sur la scène parisienne [Julia 1961, p. VI] et, d'ailleurs, presque « Allemand » [Klein 1928] ([62]). Qu'est ce qu'un « grand algébriste » en 1870 ? L'hommage rendu par Emile Picard au moment du décès de Jordan donne quelques mots clés : « Galois », « groupes », « équations », « abstraction », « généralité ». On pourrait ajouter le terme « géniale méthode de Galois » utilisé par Henri Lebesgue pour caractériser la charnière dans l'histoire de l'algèbre que représente le passage d'une science des équations à une étude « abstraite » des « groupes » ([63]) :

> Dans ses recherches, Jordan utilise la géniale méthode de Galois, dont le point essentiel est l'introduction d'un certain nombre de substitutions, déjà aperçu par Lagrange, que l'on peut attacher à chaque équation algébrique et dans lequel les propriétés des équations se reflètent fidèlement. Mais pour savoir observer dans ce miroir, il faut avoir appris à distinguer les diverses qualités des groupes de substitutions et à raisonner sur elles. C'est ce qu'à fait Jordan avec une habile ténacité et un rare bonheur ; dans son *Traité des Substitutions et des Equations algébriques*, où il a réuni et coordonné ses recherches, les propriétés des équations dérivent tout de suite de celles des groupes de substitutions.
>
> Les principales qualités des groupes qui servent à Jordan sont caractérisées par les qualités transitif ou intransitif, primitif ou imprimitif, simple ou composé. Le théorème de Jordan sur la composition des groupes est le plus connu de tous ses résultats, il entraîne cette conséquence fondamentale : il n'y a pas lieu de choisir entre les différents procédés de résolution algébrique d'une équation ; ils sont tous équivalents et conduisent aux mêmes calculs, à l'ordre près. [Lebesgue 1923, p. XX].

[59] Voir [*Comptes rendus* 1873, t. 76, p. 1302] et, dans la correspondance de Jordan, les félicitations de Cremona le 19/12/1869 [Archives Polytechnique VI2aX1855, lettre n°10], de Borchardt le 9/4/1870 [Archives Polytechnique VI2aX1855, lettre n°12], de Puiseux [Archives Polytechnique VI2aX1855, lettre n°14], de Clebsch [Archives Polytechnique VI2aX1855, lettre n°16] et de Cayley [Archives Polytechnique VI2aX1855, lettre n°26].

[60] La traduction sera attribuée à Laugel (traducteur de Riemann en 1898). Cette information est extraite d'une lettre non datée envoyée par Hermite à Jordan. Une autre origine des relations tendues entre Jordan et Hermite est probablement à attribuer au soutien apporté par Bertrand à Jordan dès la candidature de 1875 à l'Académie. Voir les critiques portées par Hermite sur Bertrand citées par [Zerner 1991].

[61] Il s'agit du théorème de Jordan–Hölder. La querelle est évitée par la médiation de Borchardt, le directeur du journal qui a vu la publication du mémoire de Netto.

Ce qui permet de voir les « vraies raisons des choses », c'est un « miroir » qui « reflète » l'étude quantitative des équations en une algèbre des « qualités », comparée à une science de la nature [Lebesgue 1923, p. XXIII]. Le *Traité* de Jordan sera présenté, dans l'histoire écrite par ses successeurs, comme matérialisant le « miroir » métaphorique symbolisant la mutation de l'algèbre par le « reflet » de deux théorèmes:

> Galois a démontré dans un mémoire célèbre [...], que chaque équation algébrique est caractérisée par un certain groupe de substitutions dans lequel se reflètent ses principales propriétés : proposition capitale qui fait dépendre la théorie toute entière des équations de celle des substitutions […].
>
> THEOREME II. – Pour qu'une équation soit résoluble par radicaux, il faut et il suffit que sa résolution se ramène à celle d'une suite d'équations abéliennes de degré premier. […] Autre énoncé du même théorème:
>
> THEOREME III. – Pour qu'une équation soit résoluble par radicaux il faut et il suffit que ses facteurs de composition soient tous premiers. [Jordan 1870, p. 385].

Les successeurs de Jordan célèbreront la synthèse du *Traité des substitutions* qui fait émerger des recherches éparses inspirées par Galois les notions et méthodes essentielles d'une théorie autonome des groupes ([64]). Par exemple, le chapitre II du *Traité* est consacré au groupe linéaire, dont la structure abstraite bénéficie d'une étude systématique (ordre, éléments générateurs, facteurs de composition etc.). Cette importance accordée au groupe linéaire ne surprendra pas le mathématicien contemporain tant elle paraît naturelle et, précisément, faire accéder la théorie des groupes « au rang de discipline autonome », c'est construire ce naturel qui fera tradition dans les manuels du XX^e siècle ([65]). Avant la synthèse de Jordan, les propriétés du groupe linéaire étaient éparpillées au sein de méthodes particulières élaborées pour la recherche des équations résolubles par radicaux. Les substitutions linéaires intervenaient notamment dans un résultat essentiel, énoncé par Galois sans démonstration :

> [Galois] a partagé les équations irréductibles en deux grandes classes : équations primitives et non primitives. Puis il a énoncé à l'égard des premières : Le degré de toute équation primitive et soluble par radicaux est une puissance d'un nombre premier. Les substitutions de son groupe sont toutes linéaires. [Jordan 1867, p. 109].

---

[62] Voir la citation de F. Klein *in* [Gispert 1991, p. 37] et plus généralement le rôle attribué à Jordan par H. Gispert dans son étude de la production mathématique française lors de la création de la SMF [Gispert 1991, p. 37].

[63] J. Dieudonné a donné un résumé dans le cadre des mathématiques qui lui étaient contemporaines des travaux algébriques de Jordan à l'occasion de la publication des œuvres de ce dernier en 1961.Voir aussi [Julia 1961, p. I-IV].

[64] La capacité de Jordan à rassembler des recherches éparses en un traité synthétique se manifeste également dans ses *Cours d'analyse* de 1893 dont H. Gispert a examiné le rôle dans le développement des fondements de l'analyse en France [Gispert, 1982].

On sait depuis Evariste Galois et Niels Henrik Abel à quelles conditions l'équation générale du $n^e$ degré est résoluble par radicaux et les travaux de Jordan concernent la détermination de toutes les équations résolubles particulières par une méthode théorique de construction des groupes résolubles maximaux du groupe des substitutions. Cette méthode, une « machinerie », une « gigantesque récurrence sur le degré $N$ de l'équation » [Dieudonné 1970, p. 168] procède d'une réduction du « genre » du groupe du général au particulier, la recherche des groupes résolubles généraux est réduite à celle de groupes particuliers définis par des « qualités » : ce sont successivement les groupes « transitifs », « primitifs » puis « linéaires ». Ces qualités correspondent aux maillons les plus simples de la chaîne de réduction. Le groupe linéaire doit son « origine » à son rôle dans cette suite de réductions et le traité de 1870 en organise les propriétés en un tout théorique et autonome. Parmi ces propriétés, l'exposé de la forme canonique des substitutions illustre le nouveau caractère théorique de la notion de groupe : d'abord élaborée en 1868 comme une méthode particulière pour les recherches sur la résolubilité des équations, la forme canonique est présentée par le *Traité* comme une réponse à une question naturelle de la théorie du groupe linéaire : « simplifier autant que possible l'expression d'une substitution » [Jordan 1870, p. 97].

Lorsque, comme nous l'avons vu dans la deuxième partie de cet article, Jordan répond en 1871-1872 aux questions posées par l'astronome Yvon-Villarceau sur la stabilité des petites oscillations, sa résolution de l'intégration des systèmes d'équations linéaires à coefficients constants participe de la même pratique de réduction d'un problème général en une suite de problèmes simples. Nécessitant la résolution d'une équation algébrique générale, la réduction canonique de Jordan ne permet pourtant pas, dans la pratique, d'intégrer un système différentiel, son application dans le champ des applications s'accompagne du transfert d'une philosophie de la généralité et d'idéaux de simplicité et d'abstraction caractéristiques des raisonnements synthétiques de la mathématique des qualités qui, pour Lebesgue caractérise la théorie des groupes [Lebesgue 1921, p. XXII].

## 3. Issues d'une querelle : la tension formes canoniques - invariants.

Deux mémoires très influents, publiés par Ferdinand Georg Frobenius en 1878 et 1879, fixeront la théorie des formes bilinéaires pour plusieurs décennies. Le premier s'ouvre par les références à Kronecker et Weierstrass, le second se clôt par une référence à Jordan. L'élaboration théorique de Frobenius se présente comme une issue à la querelle de 1874, elle doit beaucoup à la communication scientifique impulsée par celle-ci et qui se prolonge au-delà des contributions de Jordan et Kronecker avec, notamment, la publication de Darboux [1874] sur les faisceaux de formes quadratiques que nous avons déjà évoqué dans la deuxième partie et qui se rattache aux travaux d'Hermite et de Sylvester des années 1850-1860 ([66]). Frobenius déclare rechercher l'achèvement de la théorie arithmétique de Kronecker et y parvient paradoxalement en développant une idée propre à Jordan. D'une part, le mémoire écrit en mai 1877 et publié dans le tome 84 du *Journal de Crelle*, reprend la structure théorique élaborée par Kronecker, d'autre part, le titre donné au mémoire, « Ueber lineare Substitutionen und bilineare Formen » reprend l'ambition qu'avait Jordan de traiter formes et substitutions au sein d'une même théorie. Par la construction d'un calcul symbolique qui manifeste la postérité d'auteurs anglais comme Sylvester, Cayley et Smith, la théorie de Frobenius permet en 1880 de concilier deux points de vue qui s'opposaient radicalement en 1874. En référence au calcul des quaternions, la transformation linéaire des formes est traitée comme une opération symbolique de multiplication applicable aux formes comme aux substitutions puisque portant sur des « systèmes de $n^2$ valeurs ». Ce calcul permet de représenter les relations arithmétiques d'équivalence des formes définies par Kronecker en 1874, la transformation de la forme $A$ par les substitutions linéaires $x_\alpha=\sum p_{\alpha\beta}X_\beta$, et $y_\alpha=\sum q_{\alpha\beta}Y_\beta$ se formulant comme un produit symbolique de trois formes $P'AQ$ où $P=\sum p_{\alpha\beta}x_\alpha y_\beta$ et $Q=\sum q_{\alpha\beta}x_\alpha y_\beta$ ([67]). Exprimés de manière symbolique, le problème de l'équivalence des faisceaux de formes $P(rI-A)Q=rI-B$ d'une part, et celui de la similitude des substitutions $P^{-1}AP=B$ d'autre part, sont susceptibles d'une même approche. Comme nous l'avons vu, la question de l'architecture de la théorie des formes bilinéaires était un des points de désaccord entre Jordan et Kronecker. Frobenius adopte l'approche de Kronecker en mettant le problème de l'équivalence des faisceaux de formes au cœur de la théorie, la similitude des formes en corollaire. Cette organisation théorique restera prédominante jusque dans les années 1920-1930, elle est indissociable du caractère primordial donné aux calculs d'invariants et du rôle secondaire des formes canoniques ([68]).

[65] Jordan est le premier à consacrer une étude théorique au groupe linéaire dont l'importance provient de problèmes de résolubilité des équations et de fonctions modulaires. Le rôle du groupe linéaire dans les méthodes et travaux de Jordan sur les substitutions, l'origine de la forme canonique dans des problèmes de recherche de groupes résolubles et la place essentielle qu'acquiert le groupe linéaire dans les mathématiques au début du XXe siècle, avec notamment la parution du traité *Linear Groups* de Dickson en 1900 sont étudiés dans [Brechenmacher 2006a].

Dans sa publication de 1879, intitulée « Theorie der linearen Formen mit ganzen Coefficienten », Frobenius caractérise les classes de formes par des invariants polynomiaux dont la détermination ne nécessite que des « opérations rationnelles » comme l'exigeait Kronecker en 1874 ([69]). Ce mémoire s'achève avec la première démonstration de l'implication du théorème de réduction canonique de Jordan par le théorème des diviseurs élémentaires. Comme chez Weierstrass en 1868, et chez Kronecker en 1874, la recherche de la « forme normale » d'un faisceau est introduite par Frobenius pour la détermination des faisceaux de formes correspondant à un système de diviseurs élémentaires donnés [Frobenius 1879, p. 541]. Frobenius démontre que la « forme élémentaire »,

$$R = (r-a)(x_1y_1+\ldots+x_\varepsilon y_\varepsilon)-(x_1y_2+\ldots+x_{\varepsilon-1}y_\varepsilon)$$

a un unique diviseur élémentaire

---

[66] Entre 1874 et 1880, Darboux et Frobenius se préoccupent tous deux des équations différentielles du problème de Pfaff pour lesquelles les deux auteurs emploient les méthodes récentes des formes bilinéaires. Les travaux de Darboux permettent à Frobenius d'avoir accès à des méthodes élaborées par Hermite et des auteurs anglais comme Sylvester et Smith, voir [Brechenmacher 2006a, p. 246-271].

[67] La notation $P'$ désigne la transposée ${}^tP$.

[68] Consulter par exemple la communication de Dickson au congrès de Toronto de 1924 qui propose de renverser l'ordre d'exposition de la théorie des formes bilinéaires en mettant, comme le voulait Jordan en 1874, la notion de forme canonique au centre [Dickson 1924, p. 361].

[69] La démonstration de Frobenius généralise un énoncé arithmétique de Smith sur les systèmes linéaires d'équations à coefficients entiers [Smith 1861]. Ce résultat, tel que le reformule Frobenius, énonce que toute forme bilinéaire $A$ à coefficients entiers et de rang $l$ est équivalente à la forme $F=f_1x_1y_1+f_1f_2x_1y_2+\ldots+f_1f_2f_3\ldots f_lx_ly_l$. Les invariants arithmétiques $f_i$ introduits par Smith s'identifient aux diviseurs élémentaires tels que définis par Kronecker comme *p.g.c.d.* de sous déterminants. Ce sont, en des termes d'aujourd'hui, les facteurs invariants d'une matrice dont les coefficients appartiennent à un anneau principal. La similitude des propriétés arithmétiques des entiers et des polynômes permet une généralisation du résultat de Smith aux formes à coefficients polynomiaux et donc aux couples de formes bilinéaires.

$(r\text{-}a)^{\varepsilon}$

égal à son déterminant car les mineurs d'ordre *ε-1* sont égaux à 1. Les diviseurs élémentaires de la forme *R+R'+...* sont « ceux de chaque forme mis en commun ». A une suite de diviseurs élémentaires donnés correspond donc une « forme normale » obtenue par des compositions de formes du type *R*. Comme Kronecker, qui avait insisté en 1874 sur la dépendance de la forme normale à la décomposition en facteurs irréductibles du déterminant caractéristique, Frobenius présente la réduction à une forme normale comme relative à la donnée d'un « corps », au sens des travaux de Richard Dedekind sur les nombres algébriques, auquel appartiennent les coefficients des formes polynomiales. Contrairement aux arguments opposés à Jordan par Kronecker, il est cependant toujours possible d'associer des formes normales à des invariants obtenus par des procédés rationnels. Dans le cas général d'un corps de base quelconque, la forme normale est définie comme composition de formes élémentaires associées aux invariants polynomiaux *φ(r)* par le déterminant :

$$\varphi(r)=\begin{vmatrix} r+a_1 & -1 & 0 & 0 & \dots & 0 \\ a_2 & r & -1 & 0 & \dots & 0 \\ a_3 & 0 & r & -1 & \dots & 0 \\ a_4 & 0 & 0 & r & \dots & 0 \\ . & . & . & . & \dots & . \\ a_\alpha & 0 & 0 & 0 & \dots & r \end{vmatrix}$$

La forme canonique de Jordan apparait alors comme un cas singulier correspondant au cas où les coefficients du corps de base sont congruents par un nombre premier *p* ([70]). Dans ce cas, « les nombres complexes de Galois sont permis », les facteurs irréductibles du polynôme caractéristique sont linéaires et la forme normale du faisceau de formes associée est la forme *rI-J* où *J* est la forme énoncée par Jordan en 1874. Les deux théorèmes, dont l'opposition avait nourri la querelle de 1874, participent d'une identité dans la structure théorique de Frobenius, il y a désormais un seul théorème, basé sur celui de Weierstrass, dont le résultat de Jordan apparaît comme un corollaire pour la donnée d'un corps de base particulier.

A l'occasion du rapport rédigé pour sa candidature de 1881 à l'Académie de Paris, Jordan reconnaitra la priorité de la classification par Kronecker et Weierstrass des faisceaux de formes. L'influence de la théorie de Frobenius fera disparaître pour les trente ans à venir l'existence d'un théorème de Jordan de réduction canonique du paysage de la théorie des formes bilinéaires, les travaux de Kronecker seront au contraire célébrés pour leurs résultats sur l'équivalence des faisceaux singuliers ([71]). Bien que résultat d'une forte communication et mêlant des idées propres à Kronecker comme à Jordan, la réponse théorique donnée par Frobenius n'épuise cependant pas les questions d'identités posées par l'opposition des pratiques de réduction canonique et de calculs d'invariants lors de la querelle de 1874. Une forme canonique n'y est envisagée que comme un déterminant statique, représentant particulier d'une classe d'équivalence et n'apparaissant qu'à l'issue du calcul d'invariants permettant de la caractériser. A l'opposé de tels calculs d'invariants qui ne manipulent pas directement les formes, le théorème de Jordan est sous-tendu par une pratique algébrique de réduction d'un problème général en une suite de problèmes simples indissociable d'une opérationnalité de la représentation donnée aux substitutions et qui, comme le montre l'extrait suivant, permet de « voir » tout à la fois les regroupements d'indices en sous groupes structurés par la décomposition polynomiale de l'équation caractéristique et l'action de la substitution sur ces groupes ([72]) :

> Les fonctions $y_0$, $y'_0$,...;...;$y_\rho$, $y'_\rho$,...;... étant toutes distinctes, peuvent être prises pour indices indépendants à la place d'un nombre égal des indices primitifs $x, x'$, ...$x^{n-1}$. Cela fait, et $m$ étant le nombre de ces fonctions, la substitution se trouvera réduite à la forme.

$$A = \left| \begin{array}{ll} y_0 & K_0 y_0 \\ y'_0 & K_0 y'_0 \\ \ldots & \ldots \\ y_1 & K_1 y_1 \\ \ldots & \ldots \\ x^m & a_1^m x^m + b_1^m x^{m+1} + \ldots + c_1^m x^{n-1} + d_1^m y_0 + e_1^m y'_0 + \ldots + f_1^m y_1 + \ldots \\ x^{m+1} & a_1^{m+1} x^m + b_1^{m+1} x^{m+1} + \ldots + c_1^{m+1} x^{n-1} + d_1^{m+1} y_0 + e_1^{m+1} y'_0 + \ldots + f_1^{m+1} y_1 + \ldots \\ \ldots. & \ldots \\ x^{n-1} & a_1^{n-1} x^m + b_1^{n-1} x^{m+1} + \ldots + c_1^{n-1} x^{n-1} + d_1^{n-1} y_0 + e_1^{n-1} y'_0 + \ldots + f_1^{n-1} y_1 + \ldots \end{array} \right|$$

[70] C'est à dire lorsque le corps de base est un corps fini, cas considéré par Jordan en 1870.

[71] Jusqu'aux travaux de Séguier de 1907 sur la théorie des matrices, le théorème de réduction canonique sera confiné au cadre de la théorie des groupes de substitutions dans lequel il avait originellement énoncé en 1870. Au sujet de la postérité de Jordan, voir [Brechenmacher, 2006a].

[Jordan 1870, p. 117].

En associant en un même tableau les indices et leurs images, Jordan représente les étapes successives de sa pratique de réduction elle-même. Montrant l'interdépendance de la décomposition de la « substitution » et de celle des « indices », la forme que prend la substitution $A$ contient en elle-même la nouvelle substitution $C$ sur laquelle il faut itérer la réduction :

$$C = \left| \begin{array}{ll} x^m & a_1^m x^m + b_1^m x^{m+1} + \ldots + c_1^m x^{n-1} \\ x^{m+1} & a_1^{m+1} x^m + b_1^{m+1} x^{m+1} + \ldots + c_1^{m+1} x^{n-1} \\ \ldots & \ldots \\ x^{n-1} & a_1^{n-1} x^m + b_1^{n-1} x^{m+1} + \ldots + c_1^{n-1} x^{n-1} \end{array} \right|$$

Les actions de la substitution $A$ sur les « lettres », regroupées progressivement en groupes par des changements de variables linéaires, se représentent comme autant de blocs s'emboîtant chacun dans le précédent jusqu'à partager « les indices de la série » en « suites distinctes $Y_0$, $Z_0$, $u_0$,... ; $Y'_0$,$Z'_0$,... ;... que $A$ altère séparément suivant une loi "simple" $K_0Y_0$, $K_0(Z_0+Y_0)$, $K_0(u_0+Z_0)$,...;... représentée par sa forme canonique ».

$$\left| \begin{array}{ll} y_0, z_0, u_0, \ldots, y'_0, \ldots & K_0 y_0, K_0(z_0 + y_0), K_0(u_0 + z_0), \ldots, K_0 y'_0 \\ y_1, z_1, u_1, \ldots, y'_1, \ldots & K_1 y_1, K_1(z_1 + y_{10}), K_1(u_1 + z_1), \ldots, K_1 y'_1 \\ \ldots\ldots\ldots & \ldots\ldots\ldots \\ v_0 & K'_0 v_0, \ldots \\ \ldots\ldots\ldots & \ldots\ldots\ldots \end{array} \right|$$

# Conclusion.

Si la controverse de 1874 se présente comme la rencontre de deux théorèmes - celui de Weierstrass de 1868 et celui de Jordan de 1870 -, elle se caractérise surtout par l'opposition de deux pratiques – réductions canoniques et calculs d'invariants -, dont elle manifeste les identités multiples. Nous avons vu le rôle de la référence à une histoire commune impliquant des auteurs comme Lagrange, Cauchy ou Hermite pour identifier une pratique algébrique partagée par Jordan et Kronecker et consistant à transformer des systèmes linéaires par des décompositions de la forme polynomiale d'une équation algébrique. En addition de cette perspective historique large qui envisage la querelle de 1874 comme l'opposition de deux fins données à une discussion longue d'un siècle, nous avons abordé la controverse sur les plans sociologique, épistémologique et, plus généralement culturel. Sur un plan sociologique, la querelle se développe à partir d'une difficulté de communication scientifique et l'opposition des deux auteurs se déclare tout d'abord en raison de leurs appartenances à des réseaux différents. La correspondance se présente alors comme un vecteur de la part tacite des pratiques de deux savants évoluant au sein de communautés distinctes. Elle permet, sur un plan privé, de désamorcer les attaques publiques à l'origine de la querelle. La controverse n'en reste pourtant pas là et évolue en une confrontation de positions épistémologiques fortes. Analysée sur le plan épistémologique, la querelle oppose deux philosophies de l'algèbre et de la capacité de celle-ci à atteindre la généralité. Kronecker oppose à la portée théorique que veut donner Jordan à sa réduction algébrique une philosophie selon laquelle le « travail algébrique » n'a qu'une « place relative », il « est effectué au service d'autres disciplines mathématiques dont il reçoit ses fins et dont dépendent ses objectifs » [Kronecker 1874a, p. 367]. Pour le savant berlinois, la théorie des formes est de nature arithmétique et ses méthodes doivent par conséquent relever d'un calcul d'invariants obtenus par des procédés effectifs de calculs de p.g.c.d. Jordan donne au contraire une organisation algébrique à la théorie qu'il structure par l'action des groupes de substitutions, selon lui une résolution « générale » n'a de sens qu'en tant qu'elle procède de la « réduction » d'un problème jusqu'à son expression ultime qualifiée de « simple ». A l'opposé de l'effectivité revendiquée par Kronecker, le critère de réduction de Jordan présente un caractère abstrait, il demande que soient extraites toutes les racines d'une équation algébrique pour l'obtention d'une « réduite », définie par une « qualité », sa « simplicité ». Jordan reproche à Kronecker le manque de généralité de sa méthode car ses réduites peuvent encore être décomposées ; Kronecker, de son côté, reproche à Jordan sa réduction « formelle » qui nécessite des extractions de racines d'équations algébriques et que l'on ne peut réaliser en général de manière effective. Malgré son caractère abstrait, l'idéal de simplicité de Jordan, provenant de méthodes élaborées pour la recherche des groupes résolubles dans les années 1860-1870, accompagne la pratique de réduction canonique dans les applications qu'en fait le savant parisien dans diverses théories. Avant qu'un cadre disciplinaire comme celui de l'algèbre linéaire des années trente du XX^e^ siècle permette d'envisager des structures sous tendant les pratiques de Jordan ou Kronecker, l'identité algébrique ou non algébrique de telles pratiques se révèle complexe tant celles-ci ne portent le plus souvent pas de signification propre mais acquièrent des représentations diverses aux seins de différentes méthodes dans des cadres théoriques variés.

[72] Les $K_i$ désignent les racines caractéristiques. On verrait aujourd'hui dans ces regroupements les deux sous espaces stables que constituent la somme d'un espace propre et son complémentaire. Les sous espaces propres correspondants à des valeurs propres distinctes sont conjugués entre eux; leur somme constitue donc un sous espace stable $U$ par l'action de $A$ et la restriction $C$ de $A$ à $U$ permet de poursuivre la réduction. S'il ne faut pas chercher, dans la pensée de Jordan, de notion de stabilité d'un espace sous l'action d'un opérateur, le recours à des notions de l'algèbre linéaire des années 1930-1940 permet de mettre en évidence l'originalité d'une pratique indissociable d'un mode de représentation des substitutions. Pour davantage de détail sur la démonstration de Jordan, consulter [Brechenmacher, 2006a, p.167-187].

Les nombreux traités sur la « théorie des matrices » publiés dans les années trente du XX^e^ siècle donnent une place centrale à un résultat dénommé « théorème de Jordan de la décomposition matricielle ». Comme en témoigne l'adoption progressive d'une même dénomination de « forme de Jordan » sur un plan international, une nouvelle identité est donnée au théorème initialement énoncé par Jordan en 1870 et rétrospectivement interprété comme un résultat « classique » [Mac Duffee, 1933]. Cette identité est fondée sur une dualité qui donne une résolution mathématique à la tension formes canoniques/ invariants qui caractérisait la controverse de 1874 ([73]):

[73] Pour l'exemple donné ici le polynôme minimal est $\lambda^8+\alpha_1\lambda^7+\ldots+\alpha_7\lambda+\alpha_8=(\lambda-\lambda 1)^2(\lambda-\lambda_2)^3(\lambda-\lambda_3)(\lambda-\lambda_4)$.

$$A=\left\|\begin{array}{cccccccc}
\lambda_1 & 1 & 0 & & & & & 0\\
0 & \lambda_1 & 0 & & & & & \\
 & & \lambda_2 & 1 & 0 & & & \\
 & & 0 & \lambda_2 & 1 & & & \\
 & & 0 & & \lambda_2 & 0 & 0 & 0\\
 & & & & 0 & \lambda_3 & 0 & 0\\
 & & & & & & \lambda_4 & 1\\
 & & & & & & & \lambda_4
\end{array}\right\|,\ B=\left\|\begin{array}{ccccccccc}
0 & 0 & \dots & & & & \dots & 0 & -\alpha_8\\
1 & 0 & \dots & & & & & \dots & -\alpha_7\\
0 & 1 & 0 & \dots & \dots & & & & -\alpha_6\\
\dots & 0 & 1 & 0 & 0 & \dots & & & \\
 & \dots & 0 & 1 & 0 & 0 & \dots & \dots & \\
 & & \dots & 0 & 1 & 0 & 0 & \dots & \\
 & & & & 0 & 1 & 0 & & -\alpha_2\\
0 & 0 & \dots & & & \dots & 0 & 1 & -\alpha_1
\end{array}\right\|$$

D'une part la matrice $A$, « forme de Jordan », est la forme la plus simple et donne la décomposition maximale qui, pour Jordan devait caractériser les formes canoniques; d'autre part la matrice $B$, qualifiée de « forme rationnelle », est obtenue par des procédés effectifs et satisfait l'exigence de Kronecker. Aux deux théorèmes opposés en 1874, au théorème unique de la théorie de Frobenius, répondent deux formes canoniques associées à un unique théorème dont l'identité se décline sous une forme mathématique : on passe d'une forme canonique à l'autre par une méthode spécifique articulant des décompositions (arithmétique ou algébrique) d'invariants polynomiaux, des pratiques combinatoires d'extractions de sous matrices, un calcul symbolique des puissances de matrices, une arithmétique des lignes et des colonnes et le point de vue vectoriel de décomposition d'un espace en sous espaces stables par l'action d'une transformation linéaire. Les différentes démonstrations du théorème de Jordan, les applications à différents problèmes comme la recherche des matrices commutant avec une matrice donnée sont autant d'exemples de l'efficacité d'une méthode de décomposition des problèmes en une suite de sous problèmes, représentée par les formes imagées des matrices et sous-matrices, que Picard désignait encore en 1910 de « dessins » et qui, dans les années trente, envahissent les textes mathématiques.

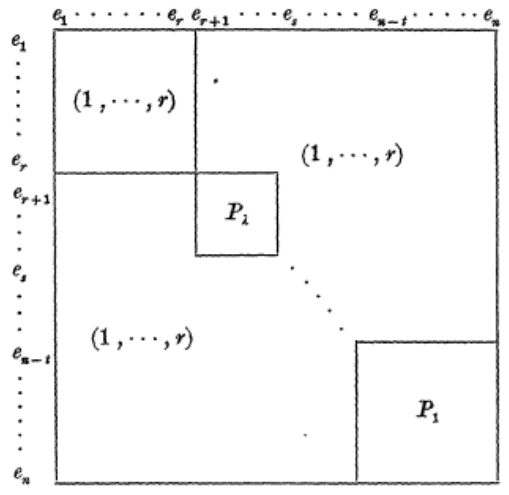

L'accent porté à cette époque sur les formes canoniques au détriment du calcul des invariants, dans le cadre d'une méthode de décomposition indissociable d'un caractère opératoire conféré à la représentation matricielle, s'accompagne de la mise en avant d'un idéal de simplicité, proche de celui qui caractérisait la position de Jordan en 1874 ([74]).

Comme nous l'avons montré dans notre thèse de doctorat, l'étude de la dynamique de la tension, mise en évidence dans cet article, entre des pratiques de réductions canoniques et de calculs d'invariants permet de suivre des fluctuations concrètes de l'élaboration mathématique entre 1870 et 1930 en prenant pour référence, non pas les structures de 1930, mais la querelle de 1874 ([75]). A la marge d'une théorie prédominante des invariants qui, ainsi que l'a montré Thomas Hawkins, donnera un caractère global à un certain nombre d'idéaux propres à l'école de Berlin, des auteurs comme Henri Poincaré, Eduard Weyr, Theodor Molien, Kurt Hensel, William Burnside, Leonard Dickson ou Léon Autonne élaborent des pratiques de réductions qui, comme celle de Jordan, s'articulent au caractère opératoire d'une représentation imagée. Ces procédés ne se limitent pas au champ aujourd'hui identifié comme correspondant à l'algèbre, l'arithmétique joue notamment un rôle important avec la combinatoire des systèmes de Kronecker ou la pratique des tableaux chez Hermite, Jordan, Poincaré et Albert Châtelet. La manière dont ces pratiques s'identifient à des réseaux - comme celui de la discussion sur « l'équation à l'aide de laquelle on détermine les inégalités séculaires des planètes » que nous avons caractérisé dans cet article - permet de mettre en évidence des aspects culturels des mathématiques comme des idéaux, des savoirs tacites, des philosophies internes, des regards historiques propres à des individus ou des communautés. En étudier les transmissions et les appropriations permet une étude fine de certains réseaux ainsi que de la manière dont les textes circulent et sont lus à des niveaux différents.

[74] Dans leur traité de 1932, Turnbull et Aitken définissent l'objet de la "théorie des matrices canoniques" comme consistant en "l'investigation systématique des types de transformations qui réduisent les matrices à leur forme la plus simple" [Aitken et Turnbull, 1932, p. 1, traduction F.B.].

[75] Un procédé similaire a été employé par C. Goldstein qui a étudié sur le long terme les relations entre analyse diophantienne et descente infinie en prenant pour référence deux théorèmes énoncés par Frenicle et Fermat. Le procédé a pour objet d'éviter l'écueil dans lequel peuvent verser les histoires à long terme en faisant « apparaître une image de la création mathématique qui ne résulte en fait que de convictions tacites : des identifications préalables peuvent ainsi venir conforter une représentation cumulative des résultats scientifiques » [Goldstein 1993, p. 180].

Questionner les identités revêtues par les pratiques algébriques permet d'enrichir l'histoire de l'algèbre linéaire. La construction d'une culture commune qu'implique la constitution de cette discipline ne saurait en effet se réduire à la problématique de l'émergence des structures. Envisager la "tâche aussi importante que difficile" consistant à approfondir l'histoire de ce que Jean Dieudonné considérait l'"un des traits marquants" de l'"algèbre moderne" et désignait sous le terme de "processus de linéarisation" [Dieudonné 1978, p. 78], nécessite une étude des procédés opératoires accompagnant les restructurations du champ des mathématiques et de l'histoire des mathématiques. Il faut notamment questionner la fabrication de l'universalité de représentations et de procédés associés à ces représentations en termes d'extension et d'unification de pratiques locales, de convergences de réseaux distinct. C'est dans ces convergences que s'élabore l'identité d'un théorème d'algèbre linéaire tel que le théorème de Jordan de la décomposition matricielle.

## ANNEXE 1. Chronologie d'une querelle.

22 décembre 1873. Note de Jordan à l'Académie de Paris [Jordan 1873].

Janvier 1874.

> 19 janvier, lecture d'un mémoire de Kronecker à l'Académie de Berlin [Kronecker 1874a].
>
> Courriers : Kronecker à Jordan, Jordan à Kronecker, Jordan à Weierstrass.

Février 1874.

> Courriers : Kronecker à Jordan, Jordan à Kronecker.
>
> 16 février, intervention non polémique de Kronecker sur les formes bilinéaires à l'Académie de Berlin [Kronecker 1874c 373-381].

Mars 1874.

> 2 mars, note de Jordan à l'Académie de Paris [Jordan 1874b].
>
> Courrier : Kronecker à Jordan (10 mars).
>
> 16 mars, intervention non polémique de Kronecker sur les formes bilinéaires à l'Académie de Berlin [Kronecker 1874b 373-381].
>
> Publication dans le *Journal de Liouville* du mémoire annoncé par la note de Jordan de 1873 [Jordan 1874a].

Avril 1874

> 23 avril, intervention de Kronecker à l'Académie de Berlin [Kronecker 1874d].
>
> 27 avril 1874, note de Kronecker à l'Académie de Paris.

Mai-juin 1874.

16 mai 1874, nouvelle intervention de Kronecker à l'académie de Berlin [Kronecker 1874c, p 382-413] et publication d'une suite (« Nachtrag ») [Kronecker 1874c] au mémoire de janvier [Kronecker 1874a].
22 juin 1874, note de Jordan à l'Académie de Paris [Jordan 1874d] suivie d'un mémoire dans le *Journal de Liouville* [Jordan 1874e].

## ANNEXE 2. La forme canonique de Jordan en *1870.*

La forme canonique des substitutions linéaires est énoncée en 1870 pour les éléments des groupes linéaires de substitutions opérant sur $p^n$ « lettres » modulo $p$ ([76]). Elle doit son origine à une suite de travaux publiés par Camille Jordan à la fin des années 1860, s'inscrivant dans le programme impulsé par les travaux de Galois et dont la synthèse réalisée par la publication en 1870 du *Traité des substitutions et des équations algébriques* sera célébrée comme l'une des origines de la théorie des groupes comme domaine mathématique propre [Julia1961, p. 1] ([77]). L'introduction du groupe linéaire au chapitre II du « Livre II Des substitutions » est exemplaire du caractère synthétique de l'ouvrage, la structure abstraite du groupe linéaire bénéficie d'une étude systématique : Jordan détermine son ordre, ses éléments générateurs, ses facteurs de composition et la détermination des classes de similitudes par la réduction d'une forme linéaire à sa forme canonique. A l'inverse de la chronologie historique, les propriétés théoriques des groupes comme le théorème de réduction canoniques seront appliquées à la détermination des équations résolubles dans le dernier livre du *Traité.*

THEOREME. *–Soit*

$$A = |x, x', \ldots\ ax+bx'+\ldots,\ a'x+b'x'+\ldots,\ \ldots|$$

*une substitution linéaire quelconque à coefficients entiers entre n indices variables chacun de 0 à p-1; Soient F, F', ... les facteurs irréductibles de la congruence de degré n*

$$\begin{vmatrix} a-K & a' & \ldots \\ b & b-K & \ldots \\ \ldots. & \ldots. & \ldots. \end{vmatrix} \equiv 0(\text{mod}.p)$$

*l, l', ...leurs degrés respectifs ; m, m', ... leurs degrés de multiplicité ;*
*On pourra remplacer les n indices indépendants x, x', ... par d'autres indices jouissant des propriétés suivantes :*
*1° Ces indices se partagent en systèmes correspondants aux divers facteurs F, F',... et contenant respectivement lm, l'm',... indices ;*
*2° Soient* $K_0, K_1, \ldots, K_{t-1}$ *les racines de la congruence irréductible F≡0 (mod.p) ; les n indices du système correspondant à F se partagent en l séries correspondantes aux racines* $K_0, K_1, \ldots, K_{l-1}$ *;*
*3° Les indices de la première série de ce système sont des fonctions linéaires des indices primitifs, dont les coefficients sont des entiers complexes formés avec l'imaginaire* $K_0$ *; ils constituent une ou plusieurs suite* $y_0, z_0, u_0, \ldots, y'_0, z'_0, \ldots; \ldots$ *(*) telles que A remplace les indices* $y_0, z_0, u_0, \ldots$ *d'une même suite respectivement par* $K_0y_0, K_0(z_0+y_0), K_0(u_0+z_0), \ldots;$
*4° Les indices de la r+1*$^{ième}$ *série sont les fonctions* $y_r, z_r, u_r, \ldots; y'_r, z'_r, \ldots; \ldots$ *respectivement conjuguées des précédente, que l'on forme en y remplaçant* $K_0$ *par* $K_r$ *; A les remplace respectivement par*

$$K_ry_r,\ K_r(z_r+y_r),\ K_r(u_r+z_r, \ldots; \ldots)$$

Cette forme simple

[76] En des termes qui nous sont contemporains, les substitutions sont définies sur un *corps* fini,
[77] Un groupe, pour Jordan, est ce que nous appellerions un groupe de permutations. Cela désigne donc, en 1870, toujours un groupe fini.

$$\left|\begin{array}{ll} y_0, z_0, u_0, \ldots, y'_0, \ldots & K_0 y_0, K_0(z_0 + y_0), K_0(u_0 + z_0), \ldots, K_0 y'_0 \\ y_1, z_1, u_1, \ldots, y'_1, \ldots & K_1 y_1, K_1(z_1 + y_{10}), K_1(u_1 + z_1), \ldots, K_1 y'_1 \\ \ldots\ldots\ldots\ldots & \ldots\ldots\ldots\ldots\ldots\ldots \\ \quad v_0 & \quad K'_0 v_0, \ldots \\ \ldots\ldots\ldots\ldots & \ldots\ldots\ldots\ldots\ldots\ldots \end{array}\right|$$

à laquelle on peut ramener la substitution $A$ par un choix d'indice convenable, sera pour nous sa forme *canonique*. [Jordan 1870, p. 127].

## ANNEXE 3. Les diviseurs élémentaires de Weierstrass.

Les couples de formes bilinéaires $(P, Q)$ sont caractérisés par des invariants définis par une décomposition en facteurs linéaires de polynômes définis par le déterminant $[P,Q]=|pP+qQ|$ et ses sous déterminants successifs. La seule décomposition en facteurs linéaires du déterminant $[P,Q]$ ne permet pas de définir un système complet d'invariants et il faut aller au-delà du déterminant en considérant les suites décroissantes d'exposants de chaque diviseur linéaire dans les « sous déterminants » successifs de $[P,Q]$. Weierstrass démontre que tout facteur commun aux sous déterminants d'ordre $n$-$x$ sera également un facteur de tous les sous déterminants d'ordre supérieur. Si on désigne par $ap+bq$ un facteur linéaire homogène d'exposant $l$ du polynôme homogène $[P,Q]$ et si $l^{(x)}$ désigne le plus grand exposant pour lequel tous les mineurs d'ordre $n$-$x$ contiennent le facteur $(ap+bq)^{l(x)}$ alors ou bien $l^{(x-1)}>l^{(x)}$ ou bien $l^{(x)}=0$. Si $r$ désigne le plus petit nombre tel que $l^{(r)}=0$, Weierstrass appelle $r$-$1$ l'index de $ap+bq$ et la suite $l, l', l'', \ldots, l^{(r-1)}$ est alors une suite strictement positive décroissante. Avec ces notations, et si l'on note $e=l-l'$, $e'=l'-l''$, $\ldots, e^{(r-1)}=l^{(r-1)}$, les relations de divisibilité entre le déterminant et ses mineurs successifs sont traduites par la décomposition suivante du facteur $(ap+bq)^l$ de $[P,Q]$ :

$$(ap+bq)^l = (ap+bq)^e (ap+bq)^{e'} \ldots (ap+bq)^{e^{(r-1)}}$$

L'ensemble des termes de cette décomposition, étendue à tous les facteurs linéaires du déterminant

$$[P,Q] = C(a_1 p + b_1 q)^{l_1} (a_2 p + b_2 q)^{l_2} \ldots (a_\rho p + b_\rho q)^{l_\rho}$$

est appelé ensemble des diviseurs élémentaires de $[P,Q]$ et noté ([78]) :

$$(a_1 p + b_1 q)^{e_1}, (a_2 p + b_2 q)^{e_2}, \ldots, (a_\rho p + b_\rho q)^{e_\rho}$$

Cette structure de décomposition du polynôme $[P,Q]$ traduit les relations de divisibilité entre le déterminant et ses sous déterminants successifs ([79]). Deux faisceaux de formes bilinéaires $pP+qQ$ et $pP'+qQ'$ sont transformables l'un en l'autre par des permutations non singulières sur les variables si et seulement si leurs ensembles de diviseurs élémentaires sont identiques ([80]).

[78] Il y a ici une ambiguïté dans la notation de Weierstrass. Les diviseurs élémentaires pour une même racine $c_1$ avaient d'abord été notés $e_1, e'_1, e''_1$ etc. Ils sont à présent notés $e_1, e_2, e_3$, etc., de sorte qu'un diviseur élémentaire $e_i$ n'est pas nécessairement associé à la racine $c_i$.

[79] Par exemple, le polynôme $P(X)=(X-3)^3(X-2)$ peut être polynôme caractéristique pour deux faisceaux de formes bilinéaires non équivalentes. Il existe trois classes d'équivalence de formes bilinéaires associées à $P$ et correspondent à trois décompositions en diviseurs élémentaires :

$$(X-3)^3(X-2)\ ;\ (X-3)^2(X-3)(X-2)\ ;\ (X-3)(X-3)(X-3)(X-2).$$

[80] Par exemple si $S = (s-c_1)^3(s-c_2)^2(s-c_3)^5$ et si les diviseurs communs aux premiers mineurs sont : $(s-c_1)(s-c_2)(s-c_3)^3$, le diviseur commun aux seconds mineurs est $(s-c_3)$, les troisièmes mineurs ont pour pgcd 1, alors les diviseurs élémentaires sont :

$$(s-c_1)^2\ ;\ (s-c_1)\ ;\ (s-c_2)\ ;\ (s-c_2)\ ;\ (s-c_3)^2\ ;\ (s-c_3)^2\ ;\ (s-c_3)\ ;\ 1$$

Avec les notations de Weierstrass :

$$e_1 = 2\ ;\ e_1'=1\ ;\ e_2 = 1\ ;\ e_2'=1\ ;\ e_3 =2;\ e_3'=2;\ e_3''=1$$

Si par les substitutions

$$(7)\begin{cases} x_1 = \sum_\gamma h_{1\gamma} u_\gamma, \ldots, x_n = \sum_\gamma h_{n\gamma} u_\gamma \\ y_1 = \sum_\gamma k_{1\gamma} v_\gamma, \ldots, y_n = \sum_\gamma k_{n\gamma} v_\gamma \end{cases}$$

où $u_1, \ldots, u_n$ et $v_1, \ldots, v_n$ sont des nouvelles variables et $h_{11}, \ldots, h_{nn}, k_{11}, \ldots, k_{nn}$ des constantes, qui ne sont soumises à aucune condition sauf que les déterminants

$$(8.) H = \begin{vmatrix} h_{11}, \ldots, h_{1n} \\ \ldots\ldots\ldots\ldots \\ h_{n1}, \ldots, h_{nn} \end{vmatrix}, K = \begin{vmatrix} k_{11}, \ldots, k_{1n} \\ \ldots\ldots\ldots\ldots \\ k_{n1}, \ldots, k_{nn} \end{vmatrix}$$

ne doivent pas être nuls, la forme $P(x_1,\ldots,x_n/y_1,\ldots,y_n)$ se transforme en une autre $P'(u_1,\ldots,u_n/v_1,\ldots,v_n)$; et de même $Q(x_1,\ldots,x_n/y_1,\ldots,y_n)$ en $Q'(u_1,\ldots,u_n/v_1,\ldots,v_n)$;alors les diviseurs élémentaires des déterminants de chaque forme $pP+qQ$, $pP'+qQ'$ coïncident.

Et inversement, quand deux couples de formes

$$P(x_1,\ldots,x_n/y_1,\ldots,y_n)\ , P'(u_1,\ldots,u_n/v_1,\ldots,v_n);$$

et

$$Q(x_1,\ldots,x_n/y_1,\ldots,y_n)\ ,\ Q'(u_1,\ldots,u_n/v_1,\ldots,v_n)$$

sont donnés dont les diviseurs élémentaires des déterminants

$$[P,Q], [P',Q']$$

coïncident ; alors il est toujours possible de déterminer les coefficients $(h_{11}, \ldots, h_{nn})$ et $(k_{11}, \ldots, k_{nn})$, telle que par les substitutions (7.) $P$ va sur $P'$ et $Q$ sur $Q'$. [Weierstrass 1868, p. 21, traduction F.B.] ([81]).

Weierstrass démontre que la forme normale simple (la forme diagonale) qui caractérise les couples de formes quadratiques peut être obtenue si et seulement si les diviseurs élémentaires sont « simples », c'est-à-dire si le degré de chaque facteur linéaire décroît d'une et une seule unité à chaque extraction de sous déterminant. Dans ce cas, la décomposition en facteur linéaire de *[P,Q]* (ou de manière équivalente, la suite des racines caractéristiques comptées avec leurs multiplicités) suffit à donner un système complet d'invariants du couple de formes. Au contraire, si le degré d'un facteur linéaire décroît de plus d'une unité à une étape de l'extraction de sous déterminants, il faut alors recourir à la suite complète des diviseurs élémentaires. Dans ce cas, il existe cependant également une « forme normale » permettant d'écrire le couple *(P,Q)* comme une juxtaposition de couples *(Φ,Ψ)* :

---

[81] Es werde durch die Substitutionen

$$(7)\begin{cases} x_1 = \sum_\gamma h_{1\gamma} u_\gamma, \ldots, x_n = \sum_\gamma h_{n\gamma} u_\gamma \\ y_1 = \sum_\gamma k_{1\gamma} v_\gamma, \ldots, y_n = \sum_\gamma k_{n\gamma} v_\gamma \end{cases}$$

wo $u_1, \ldots, u_n$ und $v_1, \ldots, v_n$ neue Veränderliche bedeuten, $h_{11}, \ldots, h_{nn}, k_{11}, \ldots, k_{nn}$ aber Constanten, welche keiner anderen Beschränkung unterworfen sind, als dass die Determinanten

$$(8.) H = \begin{vmatrix} h_{11}, \ldots, h_{1n} \\ \ldots\ldots\ldots\ldots \\ h_{n1}, \ldots, h_{nn} \end{vmatrix}, K = \begin{vmatrix} k_{11}, \ldots, k_{1n} \\ \ldots\ldots\ldots\ldots \\ k_{n1}, \ldots, k_{nn} \end{vmatrix}$$

nicht gleich Null sein dürfen, die Form

$P(x_1,\ldots,x_n/y_1,\ldots,y_n)$ in eine $P'(u_1,\ldots,u_n/v_1,\ldots,v_n)$;

und zugleich

$Q(x_1,\ldots,x_n/y_1,\ldots,y_n)$ in eine $Q'(u_1,\ldots,u_n/v_1,\ldots,v_n)$

verwandelt ; so stimmen die Determinanten der beiden Formen

$pP+qQ$, $pP'+qQ'$

in ihren Elementar-Theilern überein.

Und umgekehrt, wenn zwei Formen-Paare $P(x_1,\ldots,x_n/y_1,\ldots,y_n)$ ,$P'(u_1,\ldots,u_n/v_1,\ldots,v_n)$;und $Q(x_1,\ldots,x_n/y_1,\ldots,y_n)$ , $Q'(u_1,\ldots,u_n/v_1,\ldots,v_n)$ gegeben sind, und es stimmen die beiden Determinanten $[P,Q]$, $[P',Q']$ in ihren Elementar-Theilern überein ; so können auch stets die Coefficienten $(h_{11}, \ldots, h_{nn})$ und $(k_{11}, \ldots, k_{nn})$ so bestimmt werden, dass durch die unter (7.) angegebenen Substitutionen $P$ in $P'$ und zugleich $Q$ in $Q'$ übergeht.

$$(38).\begin{cases}\Phi = \sum_{\lambda} (X_\lambda Y_\lambda)_{e_\lambda} \\ \Psi = \sum_{\lambda} c_\lambda (X_\lambda Y_\lambda)_{e_\lambda} + \sum_{\lambda} (X_\lambda Y_\lambda)_{e_\lambda - 1}\end{cases}$$

C'est cette dernière formule qui sera à l'origine de la controverse entre Kronecker et Jordan lorsque ce dernier la déduira de sa forme canonique des substitutions linéaires.

## Annexe 4. Représentation graphique simplifiée du corpus formant la discussion sur l'équation à l'aide de laquelle on détermine les inégalités séculaires des planètes.

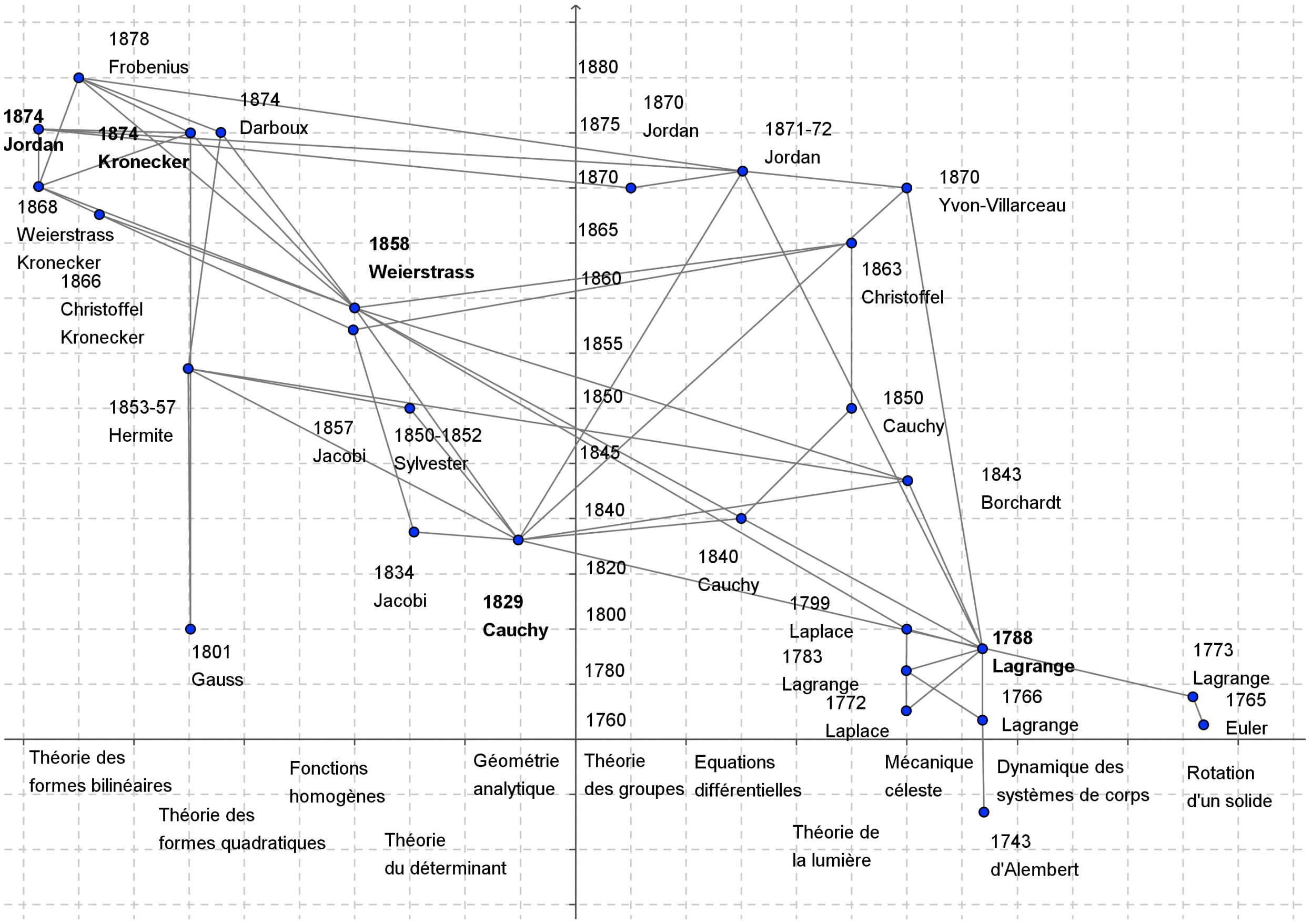

## Bibliographie.